\documentclass[12pt]{amsart}

\usepackage{amsgen,amsmath,amstext,amsbsy,amsopn,amsfonts,amssymb}
\usepackage{amsthm}
\usepackage[dvips]{graphicx}
\usepackage{epsfig}
\usepackage{ifthen}
\usepackage{amscd}
\usepackage[all]{xy}

\newfont{\fnt}{bbold10 scaled 1250}

\newcommand{\defeq}{\overset{\text{def}}{=}}

\newcommand{\Oo}{\mathcal O}

\newcommand{\Ee}{\mathcal E}

\newcommand{\Mm}{\mathfrak M}

\newcommand{\Cc}{\mathbb C}

\newcommand{\Zz}{\mathbb Z} 

\newcommand{\Id}{\mbox{\fnt 1}}

\newcommand{\CP}{\mathbb{CP}}

\newtheorem{teor}{Theorem}[section]

\newtheorem{coro}[teor]{Corollary}
\newtheorem{lema}[teor]{Lemma}
\newtheorem{prop}[teor]{Proposition}

\title[Rank stable instantons]{Rank stable instantons
over positive definite four manifolds}

\author{Jo\~ao Paulo Santos}

\begin{document}

\bibliographystyle{amsalpha}

\setlength{\parindent}{0em}
\setlength{\parskip}{1.5ex}
\newdir{ >}{{}*!/-5pt/\dir{>}}
\newdir{ |}{{}*!/-5pt/\dir{|}}

\begin{abstract}
We study the moduli space of rank stable based instantons over a connected 
sum of $q$ copies of $\CP^2$. For $c_2=1$ 
we give the homotopy type of the moduli space. For $c_2=2$ we compute the
cohomology of the moduli space.
\end{abstract}

\maketitle

\section{Introduction}

In this paper we study moduli spaces of based
$SU(r)$ instantons over a four 
manifold $X$ in the limit
when $r\to\infty$.
Interest in this rank stable limit goes back to the work of 't Hooft
\cite{tHo74}. The homotopy type of this
space was computed
in \cite{Kir94}, \cite{San95} for $X=S^4$, and
in \cite{BrSa97} for $X=\CP^2$. We have
\begin{equation}\label{eq0BS}
\Mm_k^\infty(S^4)\simeq BU(k)\ ,\ \Mm_k^\infty(\CP^2)\simeq BU(k)\times 
BU(k)
\end{equation}
where $\Mm_k$ denotes the moduli space of charge $k$ instantons.
The proofs are based on monad descriptions of the moduli spaces over
$S^4$ and $\CP^2$ (see \cite{ADHM78}, \cite{Don84},
\cite{Buc86}, \cite{Kin89}). 

In this paper we study the case where $X$ is a connected sum of $q$
copies of $\CP^2$. In \cite{Buc93}, \cite{Mat00}, it was shown that this
moduli space is isomorphic as a real analytic space to the moduli space
of holomorphic bundles on a blow up of $\CP^2$ at $q$ points, framed
at a line $L_\infty\subset\CP^2$. Under this correspondence instantons
over $S^4$ are related to holomorphic bundles on $\CP^2$ and instantons
on $\CP^2$ are related to holomorphic bundles on the blowup of 
$\CP^2$ at one point.

The results on this paper can be found in \cite{San01}. For a different
approach see \cite{KuMa99}, \cite{Buc02}
where monad descriptions for these moduli spaces were introduced.

\subsection{Results}

Fix a line $L_\infty\subset\CP^2$ and let
$x_1,\ldots,x_q\in\CP^2\setminus L_\infty$. Let $X_q$ denote the 
complex surface obtained by blowing up
$\CP^2$ at $x_1,\ldots,x_q$. 
Let $\Mm_k^r(X_q)$ be the moduli space of equivalence classes of
pairs $(\Ee,\phi)$, where $\Ee$ is
a holomorphic rank $r$ bundle over $X_q$ with
$c_1=0$ and $c_2=k$, holomorphically trivial 
at $L_\infty$, and $\phi:\Ee|_{L_\infty}\to\Oo_{L_\infty}^r$
is a holomorphic trivialization. 

For a general complex
surface $X$ the moduli space $\Mm_k^r(X)$ was defined in
\cite{Leh93}, \cite{HuLe95b},
\cite{Lub93}.

When $r_2>r_1$, there is a map $\Mm_k^{r_1}(X_q)\to\Mm_k^{r_2}(X_q)$
induced by taking direct sum with a trivial bundle: 
$\Ee\mapsto \Ee\oplus\Oo_{X_q}^{r_2-r_1}$. 
We define the rank stable moduli space
as the direct limit $\displaystyle\Mm^\infty_k(X_q)\defeq
\lim_{\substack{\longrightarrow\\r}}\Mm_k^r(X_q)$.
In this paper we study the special cases $k=1,2$.
For $k=1$ we obtain the homotopy type of the rank stable moduli space:
\begin{teor}\label{theo01}
There is a homotopy equivalence
\[
\Mm_1^\infty(X_q)\simeq BU(1)\times\left(\bigvee_{i=1}^qBU(1)\right)
\]
\end{teor}
Together with the results of \cite{Buc93} and \cite{Mat00}, this shows that
for a large class of metrics conjecture 1.1 in \cite{BrSa97} is false.

For $k=2$ we obtain the module structure of the
integer cohomology of the rank stable moduli space:

\begin{teor}\label{theo02}
Let $K_C\subset 
\Zz[x_1,x_2,x_3,x_4]\cong H^*\left(\,BU(1)^{\times4}\,\right)$ 
be the ideal generated by the product
$x_1x_2$ and let
$K_A\subset\Zz[a_1,k_1,a_2,k_2]\simeq H^*\left(\,BU(2)^{\times2}\,\right)$
be the ideal generated by $k_1,k_2$. Then,
as graded modules over $\Zz$, we have an isomorphism
\[
H^*\left(\,\Mm_2(X_q)\,\right)\cong 
\Zz[a_1,a_2]\oplus K_A^{\oplus q}\oplus K_C^{\oplus\frac{q(q-1)}2}
\]
\end{teor}

Our strategy is to  analyze the effect of the blowup on the
topology of the moduli space. This way we can relate the moduli over
$X_q$ to the moduli over $X_1$ and $X_0$, whose topology
is known (equation (\ref{eq0BS})). The plan of this paper is as follows:

In section 2 we show how the study of the moduli space $\Mm_k^r(X_q)$
can be reduced to the case where $q\leq k$.
In section 3 we recall the monad constructions of the moduli spaces
for $q=0,1$.
In section 4 we prove theorem \ref{theo01}.
Sections 5, 6 and 7 contain the proof of theorem \ref{theo02}.
The proof is based on the construction of an open cover of the moduli
space, which is carried out in section 5 and studied in detail in
sections 5 and 6.
In section 7 we use the spectral sequence associated with the open cover
(see \cite{Seg68}) to prove theorem \ref{theo02}.

\section{An open cover of $\Mm_k^r(X_q)$}

In this section we reduce the study of the moduli space $\Mm_k^r(X_q)$
to the case where $q\leq k$.

Let $I=(i_1,\ldots,i_l)$ be a multi-index and write $|I|$ for the 
number of indices.
Let $\pi_I:X_q\to X_{|I|}$ be the blow up at points
$x_j$, $j\notin I$. $\pi_I$ induces a map
\begin{equation}\label{eq2pi}
\pi_I^*:\Mm_k^r(X_{|I|})\to\Mm_k^r(X_q)
\end{equation}

The objective of this section is to prove

\begin{teor}\label{theo21}
$\left\{\,\pi_I^*\Mm_k^r(X_{|I|})\,\right\}_{|I|=k}$ 
is an open cover of $\Mm_k^r(X_q)$. Furthermore
\[
\pi_I^*\Mm_k^r(X_{|I|})\cap
\pi_J^*\Mm_k^r(X_{|J|})=
\pi_{I\cap J}^*\Mm_k^r(X_{|I\cap J|})
\]
and we have isomorphisms
\[
\xymatrix{
\Mm_k^r(X_{|I|})\ar[r]<1ex>^-{\pi_I^*}&
\pi_I^*\Mm_k^r(X_{|I|})\ar[l]<1ex>^-{\pi_{I*}}_-{\cong}
}\]
\end{teor}

From this open cover we can build a spectral sequence converging to
$H^*(\Mm_k^r(X_q))$. The case $k=2$ will be treated in section 7.
For the general case see \cite{San01}, section 4.3.

We turn now to the proof of theorem \ref{theo21}. We
begin by proving the last statement:

\begin{prop}\label{prop200}
We have isomorphisms
\[
\xymatrix{
\Mm_k^r(X_{|I|})\ar[r]<1ex>^-{\pi_I^*}&
\pi_I^*\Mm_k^r(X_{|I|})\ar[l]<1ex>^-{\pi_{I*}}_-{\cong}
}\]
where $\pi_I^*$ and $\pi_{I*}$ are inverses of each other.
We also have
\[
\pi_I^*\Mm_k^r(X_{|I|})=\left\{\Ee\in\Mm_k^r(X_q)\,|\,\Ee|_{L_i}
\mbox{ is trivial for }i\notin I\right\}
\]
\end{prop}
\begin{proof}
From theorem 3.2 in \cite{Gas97}
it follows that, 
if a bundle is trivial on the exceptional divisor
then it is also trivial on a neighborhood of the exceptional divisor.
Hence, a bundle $\Ee\to \tilde X$
on a blow up $\pi:\tilde X\to X$ is trivial on the exceptional divisor
if and only if $\tilde\Ee=\pi^*\pi_*\tilde\Ee$.
The statement of the proposition follows.
\end{proof}

\begin{proof}[Proof of theorem \ref{theo21}]
From proposition \ref{prop200}
it follows that 
\[
\pi_I^*\Mm_k^r(X_{|I|})\cap
\pi_J^*\Mm_k^r(X_{|J|})=
\pi_{I\cap J}^*\Mm_k^r(X_{|I\cap J|})
\]
To show that
$\Mm_k^r(X_q)\subset\bigcup_{|I|=k}\pi_I^*\Mm_k^r(X_k)$ we only 
need to show that:

\noindent {\bf Claim:} 
Let $\Ee\in\Mm_k^r(X_q)$, $q>k$. Then $\Ee$ is trivial in at least
$q-k$ exceptional lines.

We prove this result by induction in $q$.
Assume $\Ee$ is not trivial in
$L_1$. Let $p:X_q\to X_{q-1}$ be the blow up at $x_1$ and 
let $\Ee'=\left(\pi_*\Ee\right)^{\vee\vee}$. Then
$c_2(\Ee')<k$ so we can apply induction.
The proof is completed by noting that we cannot have bundles with
negative $c_2$ by Bogomolov inequality for framed bundles
(see \cite{Leh93}).

Finally we have to show that $\pi_I^*\Mm_k^r(X_{|I|})$ is open.
Let $H$ be an ample divisor.
Choose $N$ such that $H^i(\Ee(NH))=0$
for all ${\Ee\in\Mm_k(\tilde X)}$.
Then choose $M$ such that
$\pi_*\Ee(NH+ML)$ is locally free.
Consider the function
\[
h^1={\rm dim\,}H^1(\Ee(NH+ML)):\Mm_k(\tilde X)\to\Zz
\]
Then, from the exact sequence
\[
0\to\Ee(NH)\to\Ee(NH+ML)\to\mathcal T\to0
\]
($\mathcal T$ has support contained in $L$) we get
\[
H^2(\Ee(NH+ML))=H^2(\mathcal T)=0
\]
Now notice that
\[
H^0(\Ee(NH+ML))\cong H^0(\pi_*\Ee(NH+ML))
\]
and, since by assumption $\pi_*\Ee(NH+ML)$ is locally free and $\pi_*H$
is ample, for $N$ large enough we get 
\[
H^i(\pi_*\Ee(NH+ML))=0\mbox{ for }i>0
\]
Hence, we get that
\[
h^1=\chi(\pi_*\Ee(NH+ML))-
\chi(\Ee(NH+ML))
\]
From Riemann-Roch theorem it follows that
\[
h^1=c_2(\Ee)-c_2(\pi_*\Ee^{\vee\vee})+f(N,M,c_1(X))
\]
where $f$ does not depend on $\Ee$. The result then follows from
the upper-semicontinuity of $h^1$ (see \cite{Har77}, chapter III, 
theorem 12.8).
\end{proof}

\section{Monads}

In this section we sketch the monad description of the spaces
$\Mm_k^r(\CP^2)$ and $\Mm_k^r(\tilde\CP^2)$. We follow
\cite{Kin89}. See also \cite{BrSa00}.

Let $L_\infty\subset\CP^2$ be a rational curve and let $L$ be the
exceptional divisor.
Choose sections $x_1,x_2,x_3$ spanning $H^0(\mathcal O(L_\infty))$
and $y_1,y_2$ spanning
\mbox{$H^0(\mathcal O(L_\infty-L))$} so that
$x_3$ vanishes on $L_\infty$ and $x_1y_1+x_2y_2$ spans the kernel
of
\[
H^0(\mathcal O(L_\infty))\otimes H^0(\mathcal O(L_\infty-L))
\longrightarrow H^0(\mathcal O(2L_\infty-L))
\]

\subsection{The moduli space over $\CP^2$}

Let $W$ be a $k$-dimensional vector space.
Let $\mathcal R$ be the space of 4-tuples $m=(a_1,a_2,b,c)$ with
$a_i\in{\rm End}(W)$, $b\in{\rm Hom}(\Cc^r,W)$, $c\in{\rm Hom}(W,\Cc^r)$,
obeying the integrability condition $[a_1,a_2]+bc=0$.
For each $m=(a_1,a_2,b,c)\in\mathcal R$ we define maps $A_m,B_m$
\[
\xymatrix{
W(-L_\infty)\ar[r]^-{A_m}&
W^{\oplus2}\oplus\Cc^n\ar[r]^-{B_m}&
W(L_\infty)
}\]
by
\[
A_m=\left[\begin{array}{c}
x_1-a_1x_3\\
x_2-a_2x_3\\
cx_3
\end{array}\right]\,,\ 
B_m=\left[\begin{array}{ccc}
-x_2+a_2x_3&x_1-a_1x_3&bx_3
\end{array}\right]
\]
Then $B_mA_m=0$. The assignement 
$m\mapsto\Ee_m={\rm Ker\,}B_m/{\rm Im\,}A_m$
induces a map $f:\mathcal R\to\overline{\Mm_k^r(\CP^2)}$.

$m$ is called non degenerate if $A_m,B_m$ have maximal rank at every point
in $\CP^2$.

\begin{teor}
$f$ induces an isomorphism betwen the quotient of the space of
non degenerate points in $\mathcal R$ by the action of $Gl(W)$:
\[
g\cdot(a_1,a_2,b,c)=(g^{-1}a_1g,g^{-1}a_2g,
g^{-1}b,cg)
\]
and the moduli space $\Mm_k^r(\CP^2)$.
\end{teor}

For a proof see \cite{Don84}, proposition 1.

\begin{teor}\label{teormonadscomplete}
The algebraic quotient 
$\mathcal R/Gl(W)$ is isomorphic to the 
Donaldson-Uhlenbeck completion of the moduli space of instantons over $S^4$.
\end{teor}

For a proof see \cite{DoKr90}, sections 3.3, 3.4, 3.4.4.

For future reference we sketch here how the map
from $\mathcal R/Gl(W)$ to the
Donaldson-Uhlenbeck completion of the moduli space of instantons
is constructed (see \cite{Kin89} for details):

Let $m=(a_1,a_2,b,c)\in\mathcal R$.
A subspace $W'\subset W$ is called $b$-special with respect to $m$ if 
\begin{equation}
a_i(W')\subset W'\ (i=1,2)\ {\rm and}\ {\rm Im}\,b\subset W'.\label{1}
\end{equation}
A subspace $W'\subset W$ is called $c$-special with respect to $m$ if
\begin{equation}
a_i(W')\subset W'\ (i=1,2)\ {\rm and}\ W'\subset{\rm Ker}\,c.\label{2}
\end{equation}
$m$ is called completely reducible if for every $W'\subset W$
which is $b$-special or $c$-special, there is a complement
$W''\subset W$ which is $c$-special or $b$-special respectively.

\begin{prop}\label{PropMonadsSpecial}
Let $m=(a_1,a_2,b,c)\in\mathcal R$.
\begin{enumerate}
\item $m$ is non degenerate if and only if the only $b$-special subspace
is $W$ and the only $c$-special subspace is $0$;
\item For every $m$, the orbit of $m$ under $Gl(W)$ contains in its closure
a canonical completely reducible orbit and completely reducible orbits
have disjoint closures;
\item If $m$ is completely reducible then, after acting with some
$g\in Gl(W)$ we can write
\[
a_i=\begin{bmatrix}a_i^{red}&0\\0&a_i^{\Delta}\end{bmatrix}\ ,\ 
b=\begin{bmatrix}b^{red}\\0\end{bmatrix}\ ,\ 
c=\begin{bmatrix}c^{red}&0\end{bmatrix}
\]
where $(a_1^{red},a_2^{red},b^{red},c^{red})$ is non-degenerate
and the matrices $a_1^\Delta,a_2^\Delta$ can be
simultaneously diagonalized. Such a configuration is equivalent
to the following data:
\begin{itemize}
\item An irreducible integrable configuration 
$(a_1^{red},a_2^{red},b^{red},c^{red})$
corresponding to a bundle with $c_2=l\leq k$;
\item $k-l$ points in $\Cc^2=\CP^2\setminus L_\infty$ 
given by the eigenvalue pairs of
$ a_1^\Delta,a_2^\Delta$
\end{itemize}
This is precisely the Donaldson-Uhlenbeck completion.
\end{enumerate}
\end{prop}

\subsection{The moduli space over $\tilde\CP^2$}

Let $\tilde{\mathcal R}$ be the space of 5-tuples 
$\tilde m=(a_1,a_2,d,b,c)$ where
$a_i\in{\rm Hom}(W,V)$, $d\in{\rm Hom}(V,W)$,
$b\in{\rm Hom}(\Cc^r,V)$, $c\in{\rm Hom}(W,\Cc^r)$,
such that $a_1(W)+a_2(W)+b(\Cc^r)=V$,
obeying the integrability condition $a_1da_2-a_2da_1+bc=0$.
For each $\tilde m=(a_1,a_2,d,b,c)\in\tilde{\mathcal R}$ 
we define maps $A_{\tilde m},B_{\tilde m}$
\begin{multline*}
W(-L_\infty)\oplus V(L-L_\infty)\overset{A_{\tilde m}}{\longrightarrow}
\left(V\oplus W\right)^{\oplus2}\oplus\Cc^n
\overset{B_{\tilde m}}{\longrightarrow}\\
\rightarrow V(L_\infty)\oplus W(L_\infty-L)
\end{multline*}
by
\[
A_{\tilde m}=\left[\begin{array}{cc}
a_1x_3&-y_2\\
x_1-da_1x_3&0\\
a_2x_3&y_1\\
x_2-da_2x_3&0\\
cx_3&0
\end{array}\right]\,,\ 
B_{\tilde m}=\left[\begin{array}{ccccc}
x_2&a_2x_3&-x_1&-a_1x_3&bx_3\\
dy_1&y_1&dy_2&y_2&0
\end{array}\right]
\]
Then $B_{\tilde m}A_{\tilde m}=0$. 
The assignement 
$\tilde m\mapsto\Ee_{\tilde m}={\rm Ker\,}B_{\tilde m}/{\rm Im\,}A_{\tilde m}$
induces a map $\tilde f:\tilde{\mathcal R}\to\overline{\Mm_k^r(\tilde\CP^2)}$.

A point $\tilde m\in\tilde{\mathcal R}$ is called non-degenerate 
if $A_{\tilde m}$ and $B_{\tilde m}$ 
have maximal rank at every point in $\tilde\CP^2$. 

\begin{teor}
The map $\tilde f$ induces an isomorphism
betwen the quotient of the space of
non degenerate points in $\tilde{\mathcal R}$ 
by the action of $Gl(V)\times Gl(W)$:
\[
(g_0,g_1)\cdot(a_1,a_2,b,c,d)=(g_0^{-1}a_1g_1,g_0^{-1}a_2g_1,
g_0^{-1}b,cg_1,g_1^{-1}dg_0)
\]
and the moduli space $\Mm_k^r(\tilde\CP^2)$.
\end{teor}

See \cite{Kin89} for a proof.

Consider the algebraic quotient
$\tilde{\mathcal R}/Gl(V)\times Gl(W)$. This space is a completion
of the moduli space $\Mm_k^r(\tilde\CP^2)$. We proceed
to give an interpretation of the points in this completion in terms
of the Donaldson-Uhlenbeck completion. See \cite{Kin89} for details.

Let $\tilde m=(a_1,a_2,d,b,c)$. 
Let $V'\subset V$ and $W'\subset W$ and assume
${\rm dim}\,V'={\rm dim}\,W'$. The pair $(V',W')$
is called $b$-special with respect to $\tilde m$ if
\begin{equation}
a_i(W')\subset V'\ (i=1,2),\ d(V')\subset W'
\ {\rm and}\ {\rm Im}\,b\subset V'\label{11}
\end{equation}
The pair $(V',W')$ is called $c$-special with respect to $\tilde m$ if
\begin{equation}
a_i(W')\subset V'\ (i=1,2),\ d(V')\subset W'
\ {\rm and}\ W'\subset{\rm Ker}\,c\label{22}
\end{equation}
$\tilde m$ is called completely reducible if for every pair
$(V',W')$ which is either $b$-special or $c$-special, there are complements
$V'',W''$ to $V'$ and $W'$ such that the pair $(V'',W'')$ is $c$-special
or $b$-special respectively.

\begin{prop}\label{prop1nnX1}
Let $\tilde m=(a_1,a_2,d,b,c)\in\tilde R$.
\begin{enumerate}
\item $\tilde m$ is non-degenerate if and only if the only 
$b$-special pair is $(V,W)$ and the only $c$-special pair is $(0,0)$;
\item For every $\tilde m$, the orbit of $\tilde m$ under $Gl(V)\times Gl(W)$
contains in its closure a canonical completely reducible orbit and
completely reducible orbits have disjoint closures;
\item If $\tilde m$ is completely reducible then, after acting with some
$(g_0,g_1)\in Gl(V)\times Gl(W)$, we can write
\[
a_i=\begin{bmatrix}a_i^{red}&0\\0&a_i^{\Delta}\end{bmatrix}\ ,\ 
d=\begin{bmatrix}d^{red}&0\\0&d^\Delta\end{bmatrix}
b=\begin{bmatrix}b^{red}\\0\end{bmatrix}\ ,\ 
c=\begin{bmatrix}c^{red}&0\end{bmatrix}
\]
where $(a_1^{red},a_2^{red},d^{red},b^{red},c^{red})$ is non-degenerate
effective and integrable and the matrices $a_1^\Delta,a_2^\Delta,d^\Delta$ 
can be simultaneously diagonalized.
Such a configuration is equivalent to the following data:
\begin{itemize}
\item An irreducible configuration
$(a_1^{red},a_2^{red},d,b^{red},c^{red})$
associated to a bundle with $c_2=l\leq k$;
\item $k-l$ points in the blow up
$\widetilde{\Cc^2}$ of $\Cc^2$ at the origin. This points are
determined as follows:
$a_1^\Delta,a_2^\Delta,d^\Delta$ determine $k-l$ unique points
$(\lambda_1^r,\lambda_2^r),[\mu_1^r,\mu_2^r]\in
\widetilde{\Cc^2}$ corresponding to vectors $v_1,\ldots,v_{k-l}$
such that $da_iv^r=\lambda_i^rv^r$
($\lambda_1,\lambda_2$ are the eigenvalue pairs of $da_1,da_2$) and
$(\mu_1^ra_1+\mu_2^ra_2)v^r=0$.
\end{itemize}
\end{enumerate}
\end{prop}

\subsection{Direct image}

In this section we gather some results concerning
the direct image map $\pi_*$ induced by the blowup
map $\pi:\tilde\CP^2\to\CP^2$.

\begin{prop}
Let $\pi_\#:\tilde{\mathcal R}\to\mathcal R$ be given by
$\pi_\#(a_1,a_2,d,b,c)=(da_1,da_2,db,c)$. Let $\tilde m\in\tilde{\mathcal R}$,
$m=\pi_\#\tilde m$. Then $\Ee_{\tilde m}|_{\tilde\CP^2\setminus L}$
is isomorphic to $\Ee_m|_{\CP^2\setminus[0,0,1]}$.
\end{prop}

For the proof see \cite{San02}, proposition 5.6.

\begin{prop}\label{prop37S0}
Let $S_0\Mm_1^r(\tilde\CP^2)=
\{(\Ee,\phi)\in\Mm_1^r(\CP^2):(\pi_*\Ee)^{\vee\vee}=\Oo_{\CP^2}^r\}$. Then
\begin{enumerate}
\item $m\in S_0\Mm_1^r(\CP^2)$ if and only if $m$ is of the form
$(a_1,a_2,0,b,c)$.
\item The inclusion $S_0\Mm_1\to\Mm_1$ is a homotopy equivalence.
\end{enumerate}
\end{prop}
\begin{proof}
First we observe that
$\Mm_1(\tilde\CP^2)=S_0\Mm_1(\tilde\CP^2)\cup\pi^*_\emptyset\Mm_1(\CP^2)$.
Now $m\in\pi^*_\emptyset\Mm_1(\CP^2)$ if and only if $d$
is an isomorphism (see \cite{Kin89}). The first statement follows.
The second statement follows easily from the first: just consider the
homotopy $(a_1,a_2,d,b,c)\mapsto(a_1,a_2,td,b,c)$.
\end{proof}

\begin{prop}\label{proptau1}
Let $x=[x_1,x_2,1]\in\CP^2$ and let $\pi_x:\tilde\CP^2\to\CP^2$ be the blow up
at $x$. Then the map $\pi_x^*:\Mm_k(\CP^2)\to\Mm_k(\tilde\CP^2)$ is given by
\[
[a_1,a_2,b,c]\mapsto[a_1-x_1\Id,a_2-x_2\Id,\Id,b,c]
\]
\end{prop}
\begin{proof}
For $x=[0,0,1]$ see \cite{BrSa00}. For the general case consider
the translation
$[w_1,w_2,w_3]\mapsto[w_1-x_1w_3,w_2-x_2w_3,w_3]$.
This induces a map $\tau:\Mm_k(X_0)\to\Mm_k(X_0)$ given by
\[
[a_1,a_2,b,c]\mapsto[a_1-x_1\Id,a_2-x_2\Id,b,c]
\]
The result follows.
\end{proof}

\section{The charge one moduli space}

The objective of this section is to prove theorem \ref{theo01}:

\begin{teor}\label{theo25}
There is a homotopy equivalence
\[
\Mm_1^\infty(X_q)\simeq BU(1)\times\left(\bigvee_{i=1}^qBU(1)\right)
\]
\end{teor}

From theorem \ref{theo21} it follows that 
\[
\Mm_1^r(X_q)=\bigcup_{l=1}^q\pi_l^*\Mm_1^r(X_1)
\]
and, for any $i\neq j$,
\[
\pi_i^*\Mm_1^r(X_1)\cap\pi_j^*\Mm_1^r(X_1)=\pi_\emptyset^*\Mm_1^r(X_0)
\]
We begin by studying the maps 
$\pi_\emptyset^*\Mm_1(X_0)\to\pi_i^*\Mm_1(X_1)$.

\begin{lema}\label{prop27}
Let $\iota_1:\CP^r\to\CP^r\times\CP^r$ be the inclusion into the first factor:
$\iota_1([u])=([u],*)$, where $*$ denotes the base point. Then
there are homotopy equivalences $h_0:\CP^\infty\to\Mm_1^\infty(X_0)$
and $h_1:\CP^\infty\times\CP^\infty\to\Mm_1^\infty(X_1)$ 
such that the following diagram
\[
\xymatrix{
\pi_\emptyset^*\Mm_1^\infty(X_0)\ar[r]&\pi_i^*\Mm_1^\infty(X_1)\\
\Mm_1^\infty(X_0)\ar[u]^{\pi_\emptyset^*}_{\cong}\ar[r]^-{\pi^*}&
\Mm_1^\infty(X_1)\ar[u]^{\pi_i^*}_{\cong}\\
\CP^\infty\ar[u]^{h_0}_{\cong}\ar[r]^-{\iota_1}&
\CP^\infty\times\CP^\infty\ar[u]^{h_1}_{\cong}
}\]
is homotopy commutative.
\end{lema}
\begin{proof}
We will use the monad descritpion of
$\Mm_1^r(X_1),\Mm_1^r(X_0)$.
We define define the following maps: 
\begin{align*}
&p_0:\Mm_1^r(X_0)\to\CP^r&&p_0:[a_1,a_2,b,c]\to[b]\\
&p_1:\Mm_1^r(X_1)\to\CP^r\times\CP^r&&p_1:[a_1,a_2,d,b,c]\to\left([b],
\left[\frac{\bar c^t}{\|c\|^2}\right]\right)\\
&\Delta:\CP^r\to\CP^r\times\CP^r&&\Delta:[u]\to([u],[u])\\
&f:\CP^\infty\times\CP^\infty\to\CP^\infty\times\CP^\infty&&
f:(x,y)\mapsto(x,xy^{-1})
\end{align*}
where to define $f$ we observe that 
$\CP^\infty=BU(1)$ is homotopic to the free abelian group on $U(1)$.
Now observe that the diagram
\[
\xymatrix{
\pi_\emptyset^*\Mm_1^\infty(X_0)\ar[r]&\pi_i^*\Mm_1^\infty(X_1)\\
\Mm_1^\infty(X_0)
\ar[u]^{\pi_\emptyset^*}_{\cong}\ar[r]^-{\pi^*}\ar[d]^-{p_0}&
\Mm_1^\infty(X_1)\ar[u]^{\pi_i^*}_{\cong}\ar[d]^-{p_1}\\
\CP^\infty\ar[r]^-{\Delta}\ar[rd]^-{\iota_1}&
\CP^\infty\times\CP^\infty\ar[d]^f\\
&\CP^\infty\times\CP^\infty
}\]
is homotopy commutative and the maps $p_0,p_1,f,\pi_\emptyset^*,\pi_i^*$ 
are homotopy equivalences. The statement of the lemma then follows by writing
$h_0=p_0^{-1}$ and $h_1=p_1^{-1}f^{-1}$,
where $p_0^{-1},p_1^{-1},f^{-1}$ are the homotopy inverses.
\end{proof}

We are ready to prove theorem \ref{theo25}.

\begin{proof}
Let $C$ be the cone on $q$ points $v_1,\ldots,v_q$. Let
\[
M=\frac{
\left(\displaystyle\coprod_{i=1}^q BU(1)\times BU(1)\times\{v_i\}\right)
\coprod\left(\,BU(1)\times C\,\right)}{
(\,[u]\,,\,v_i\,)\sim(\,\iota_1([u])\,,\,v_i)}
\]
\begin{enumerate}
\item First we show that $M$ is homotopically equivalent to $\Mm_1(X_q)$.
Denote the points in $C$ by
\[
[t,v_i]\in C=\frac{[0,1]\times\coprod_i\{v_i\}}{(0,v_i)\sim(0,v_j)\sim *}
\]
Define a map
\[
\zeta:\left(\coprod_{i=1}^qBU(1)\times BU(1)\times\{v_i\}\right)
\coprod BU(1)\times C\to \Mm_1(X_q)
\]
as follows:
\begin{gather*}
BU(1)\times BU(1)\times\{v_i\}\ni
([u_1],[u_2],v_i)\mapsto\pi_i^*h_1([u_1],[u_2])\\
BU(1)\times C\ni
([u],[t,v_i])\mapsto\pi_\emptyset^*h_0([u])
\mbox{ for }t<\frac13\\
BU(1)\times C\ni
([u],[t,v_i])\mapsto\pi_i^*h_1\iota_1([u])
\mbox{ for }t>\frac23
\end{gather*}
For $\frac13\leq t\leq\frac23$ use the homotopy between 
$\pi_\emptyset^*h_0$ and $p_i^*h_1\iota_1$ from lemma \ref{prop27}.

$\zeta$ descends to the quotient to give a map $\zeta:M\to\Mm_1^\infty(X_q)$.
We want to apply Whitehead theorem
to show $\zeta$ is a homotopy equivalence.
The van Kampen theorem implies both $M$ and $\Mm_1(X_q)$
are simply connected hence we only have to show $\zeta$ is an isomorphism
in homology groups. We prove it by induction in $q'=1,\ldots,q$.

We apply the five lemma
to the Meyer-Vietoris long exact sequence corresponding to open 
neighborhoods of the sets
\begin{align*}
&\pi_{q'+1}^*\Mm_1^\infty(X_1),\,
\pi_{(1,\ldots,q')}^*\Mm_1^\infty(X_{q'})\subset\Mm_1^\infty(X_q)\\
&BU(1)\times BU(1)\times\{v_i\},\,
\bigcup_{l=1}^{q'}BU(1)\times BU(1)\times\{v_l\}\subset M
\end{align*}
using the fact that the restrictions
\begin{align*}
&\zeta:BU(1)\times BU(1)\times\{v_i\}\to\pi_i^*\Mm_1^\infty(X_1)\\
&\zeta:BU(1)\times C\to\pi_\emptyset^*\Mm_1^\infty(X_0)
\end{align*}
are homotopy equivalences.
It follows that $\zeta$ induces isomorphisms in all homology groups.
\item To conclude the proof we only have to show that
$M$ is homotopically equivalent to 
\[
BU(1)\times\left(\bigvee_{i=1}^qBU(1)\right)=
\frac{\displaystyle\coprod_{l=1}^qBU(1)\times BU(1)\times\{v_l\}}{
(x,*,v_i)\sim(x,*,v_j)}
\]
where $*\in BU(1)$ is the base point.
Define an open cover of $BU(1)\times\left(\bigvee_iBU(1)\right)$ by
$U_i=BU(1)\times BU(1)\times\{v_i\}$. 
Then the claim
is a special case of proposition 4.1 in \cite{Seg68}.
The result can also pe proved as in step 1.
\end{enumerate}
\end{proof}

\section{An open cover of $\Mm_2^\infty(X_2)$}

The objective of this section is to describe an open cover of
$\Mm_2^\infty(X_2)$.
We will adopt, in this section and the next, the following notation:
Denote the blow up points by $x_L,x_R\in X_0$.
Let $\pi:X_2\to X_0$ be the blow up map at
$x_L,x_R$. By abuse of notation we will denote by
$\pi_L$ the maps $X_2\to X_1$ and $X_1\to X_0$ corresponding to the blow
up at $x_L$ and in the same way $\pi_R$ will denote the blow up at $x_R$.
We have the diagram
\[
\xymatrix{
&X_2\ar[dl]_-{\pi_L}\ar[dr]^-{\pi_R}\\
X_{1R}\ar[dr]_-{\pi_R}&&X_{1L}\ar[dl]^-{\pi_L}\\
&X_0
}\]
of blow up maps where $X_{1L}\cong X_{1R}\cong X_1$. Denote by $L_L$ and
$L_R$ the exceptional divisors above $x_L$ and $x_R$ respectively.
Again, by abuse of notation we identify $L_L\subset X_2$ with
$L_L\subset X_{1L}$ and the same for $L_R$.

Write $x_L=[x_{1L},x_{2L},1],\,x_R=[x_{1R},x_{2R},1]$,
$x_L,x_R\in X_0=\CP^2$. 
Since $x_L\neq x_R$ we may assume without loss
of generality that $x_{1L}\neq x_{1R}$.

Let $z_i=x_{iR}-x_{iL}$.
$z_1,z_2$ determine a point 
$([z_1,z_2,1],[z_1,z_2])\in X_1\setminus L_\infty=
\widetilde{\CP^2}\setminus L_\infty\subset\CP^2\times\CP^1$.

We are ready to state the main theorem of this section:

\begin{teor}\label{teor51opencover}
Let
\begin{align*}
&A_L=\pi_R^*\Mm_2(X_{1L})=
\left\{\Ee\in\Mm_2(X_2)\,:\,\Ee|_{L_R}\mbox{ is trivial}\right\}\\
&A_R=\pi_L^*\Mm_2(X_{1R})=
\left\{\Ee\in\Mm_2(X_2)\,:\,\Ee|_{L_L}\mbox{ is trivial}\right\}
\end{align*}
and let $C=\Mm_2(X_2)\setminus(A_L\cup A_R)$.

Let $N_L\subset\Mm_2(X_{1L})$ be the set of 
non-degenerate configurations $m=(a_1,a_2,d,b,c)$ such that 
the eigenvalues of $da_1$ (equal to the eigenvalues of $a_1d$)
are in a $\delta$ neighborhood of $0,z_1$. 
In a similar way define $N_R\subset\Mm_2(X_{1R})$.
Let $N_2=\pi_R^*N_L\cup\pi_L^*N_R\cup C$.

Then $\{A_L,A_R,N_2\}$ is an open cover of $\Mm_2^\infty(X_2)$.
There are homotopy equivalences
\begin{enumerate}
\item $A_L\simeq A_R\simeq BU(2)\times BU(2)$
\item $C\simeq BU(1)\times BU(1)\times BU(1)\times BU(1)$
\item $A_L\cap A_R\simeq BU(2)$
\item $A_L\cap N_2\simeq N_L\simeq A_R\cap N_2\simeq N_R\simeq 
BU(1)\times BU(1)\times BU(1)$
\item $A_L\cap A_R\cap N_2\simeq BU(1)\times BU(1)$
\item $N_2\simeq C$
\end{enumerate}
\end{teor}

From this open cover we get, in a standard way (see \cite{Seg68}),
a spectral sequence:

\begin{coro}\label{corospseq2}
There is a spectral sequence
converging to the cohomology of $\Mm_2^r(X_2)$ with $E_1$ term
\begin{align*}
&E_1^{0,n}=H^n(A_L)\oplus H^n(A_R)\oplus H^n(N_2)\\
&E_1^{1,n}=H^n(A_L\cap A_R)\oplus H^n(A_L\cap N_2)\oplus H^n(A_R\cap N_2)\\
&E_1^{2,n}=H^n(A_L\cap A_R\cap N_2)
\end{align*}
\end{coro}
In the next section we will study the $d_1$ differential of this
spectral sequence.

We turn now to the proof of theorem \ref{teor51opencover}.
We will delay the proof that $N_2$ is open and
begin by proving the homotopy equivalences (1), (2) and (3):

\begin{prop}\label{prop50}
$A_L,A_R$ are open sets,
\[
C=\left\{[\Ee,\phi]\in\Mm_2(X_2)\,:\,c_2\left(\,
(\pi_{i*}\Ee)^{\vee\vee}\,\right)=1\,,\,i=L,R\right\}
\]
and the following maps are isomorphisms:
\begin{align*}
&\pi_R^*:\Mm_2(X_{1L})\to A_L\subset\Mm_2(X_2)\\
&\pi_L^*:\Mm_2(X_{1R})\to A_R\subset\Mm_2(X_2)\\
&\pi_{R*}^{\vee\vee}\times\pi_{L*}^{\vee\vee}:
C\to S_0\Mm_1(X_{1L})\times S_0\Mm_1(X_{1R})\\
&\pi^*:\Mm_2(X_0)\to A_L\cap A_R\subset\Mm_2(X_2)
\end{align*}
where $\pi_{i*}^{\vee\vee}(\Ee)\defeq(\pi_{i*}\Ee)^{\vee\vee}$.
\end{prop}
\begin{proof}
The isomorphisms for $A_L,A_R,A_L\cap A_R$ follows from
theorem \ref{theo21}. That theorem also implies $A_L,A_R$ are open.

It remains to look at the map
$\pi_{R*}^{\vee\vee}\times\pi_{L*}^{\vee\vee}:
C\to S_0\Mm_1(X_{1L})\times S_0\Mm_1(X_{1R})$. 
The continuity of this map was proved in proposition
3.1 in \cite{San02}.
We will construct an inverse for 
$\pi_{R*}^{\vee\vee}\times\pi_{L*}^{\vee\vee}$.
Let $(\Ee_L,\phi_L)\in S_0\Mm_1(X_{1L})$, $(\Ee_R,\phi_R)\in S_0\Mm_1(X_{1R})$.
Hartog's theorem implies there are unique extensions of
$\phi_L,\phi_R$ to maps
\[
\phi_L:\Ee_L|_{X_0\setminus\{x_L\}}\to\Oo_{X_0\setminus\{x_L\}}^r\ ,\ 
\phi_R:\Ee_R|_{X_0\setminus\{x_R\}}\to\Oo_{X_0\setminus\{x_R\}}^r
\]
These maps induce an isomorphism $\Ee_L\cong\Ee_R$ over
$X_0\setminus\{x_L,x_R\}$ which we use to glue $\Ee_L,\Ee_R$ and
obtain a bundle $\Ee\to X_2$. The continuity of this map was proved
in proposition 3.3 in \cite{San02}. This concludes the proof.
\end{proof}

We observe the following identity:

\begin{prop}\label{proptau2}
Let $\tau:\Mm_k(X_0)\to\Mm_k(X_0)$ be defined by 
\[
\tau(a_1,a_2,b,c)=(a_1-x_{1L}\Id,a_2-x_{2L}\Id,b,c)
\]
Let $m_1,m_2\in\Mm_1(X_0)$. Then
$\pi_L^*(m_1\boxplus_0 m_2)=\pi_L^*m_1\boxplus_L\tau(m_2)$.
\end{prop}
\begin{proof}
It follows easily from proposition \ref{proptau1}.
\end{proof}

Before we continue we need the lemma

\begin{lema}\label{lemm511}
Let $m\in\overline{\Mm_2(X_1)}$ and let $(a_1,a_2,d,b,c)$ be 
the configuration associated to $m$. 
The following are equivalent:
\begin{enumerate}
\item $\Ee_m$ is in the image of $\pi_{R*}:C\to\overline{\Mm_2(X_1)}$;
\item $cdb=0$ and
the eigenvalues of $da_i$ (equal to the ones of $a_id$) are
$0$ and $z_i$;
\item After a change of basis we can write
\[
a_1=\begin{bmatrix}a_1'&0\\0&z_1\end{bmatrix}\ ,\ 
a_2=\begin{bmatrix}a_2'&\frac{b'c''}{z_1}\\
-\frac{b''c'}{z_1}&z_2\end{bmatrix}\ ,\ 
d=\begin{bmatrix}0&0\\0&1\end{bmatrix}\ ,\ 
b=\begin{bmatrix}b'\\b''\end{bmatrix}\ ,\ 
c=\begin{bmatrix}c'&c''\end{bmatrix}
\]
with $c''b''=0$. 
\end{enumerate}
\end{lema}
\begin{proof}
We will show that $1\Rightarrow 2$, $2\Rightarrow 3$ and $3\Rightarrow 1$.
\begin{enumerate}
\item[$(1\Rightarrow 2)$]
Suppose $\Ee_m=\pi_{R*}\tilde \Ee$, $\tilde \Ee\in C$. 
Then, by proposition \ref{prop50}, 
$\Ee_m^{\vee\vee}\in S_0\Mm_1(X_{1L})$
and $\Ee$ is not locally free at the blow up point $x_R$.
So, from proposition \ref{prop37S0},
$\Ee^{\vee\vee}$ corresponds to a configuration of the form
$[a_1',a_2',0,b',c']$.

Since $m$ is degenerate, by proposition \ref{prop1nnX1}
after a change of basis 
it can be written in one of two forms, corresponding to the two
types of special pairs:
\begin{enumerate}
\item[1 ($b$-special):]
\[
a_i=\begin{bmatrix}a_i'&*\\0&a_i''\end{bmatrix}\ ,\ 
d=\begin{bmatrix}d'&*\\0&d''\end{bmatrix}\ ,\ 
b=\begin{bmatrix}b'\\0\end{bmatrix}\ ,\ 
c=\begin{bmatrix}c'&c''\end{bmatrix}
\]
in which case the configuration is equivalent to the 
completely reducible configuration (see proposition \ref{prop1nnX1})
\[
(a_1',a_2',d',b',c')\oplus(a_1'',a_2'',d'',0,c'')
\]
corresponding to an ideal bundle with singularity at
$(a_1''d'',a_2''d'')$ and charge one bundle given by $(a_1',a_2',d',b',c')$.
So we should have $d'=0$ and $a_i''d''=z_i$.
\item[2 ($c$-special):]
\[
a_i=\begin{bmatrix}a_i'&0\\ *&a_i''\end{bmatrix}\ ,\ 
d=\begin{bmatrix}d'&0\\ *&d''\end{bmatrix}\ ,\ 
b=\begin{bmatrix}b'\\b''\end{bmatrix}\ ,\ 
c=\begin{bmatrix}c'&0\end{bmatrix}
\]
in which case the configuration is equivalent to the
completely reducible configuration
\[
(a_1',a_2',d',b',c')\oplus(a_1'',a_2'',d'',b'',0)
\]
corresponding to an ideal bundle with singularity at
$(a_1''d'',a_2''d'')$ and charge one bundle given by $(a_1',a_2',d',b',c')$.
So we should have $d'=0$ and $a_i''d''=z_i$.
\end{enumerate}
In both cases the eigenvalues of $da_i$ are $0,z_i$ and
$cdb=0$.
\item[$(2\Rightarrow 3)$]
Now assume the configuration
$(a_1,a_2,d,b,c)$ satisfies 2. Fix a basis of eigenvectors $v_0,v_1\in V$
of $a_1d$ and $w_0,w_1\in W$ of $da_1$ with $v_0,w_0$ corresponding
to the eigenvalue $0$. Normalize $v_1,w_1$ so that $dv_1=w_1$. Then
\begin{gather}\label{eq5nboxplusL}
a_1=\begin{bmatrix}
a_1'&0\\0&a_1''\end{bmatrix}\ ,\ 
a_2=\begin{bmatrix}
a_2'&\frac{b'c''}{a_1''-d'a_1'}\\\frac{b''c'}{d'a_1'-a_1''}&a_2''
\end{bmatrix}\\
d=\begin{bmatrix}
d'&0\\0&1\end{bmatrix}\ ,\ 
b=\begin{bmatrix}b'\\b''\end{bmatrix}\ ,\ 
c=\begin{bmatrix}c'&c''\end{bmatrix}
\end{gather}
From $cdb=0$ we get $(b'c'')(b''c')=0$.
If $b'c''=0$ then $a_2$ is lower triangular. If
$b''c'=0$ then $a_2$ is upper triangular. In both cases
the diagonal entries of $a_2d$ are its eigenvalues.
Hence, the condition about the eigenvalues of $a_1d$
and $a_2d$ yields the equations
\[
a_1'd'=a_2'd'=0\ ,\ a_1''=z_1\ ,\ a_2''=z_2
\]
Since $a_1(W)+a_2(W)+b(\Cc^r)=V$ we must have $d'=0$.
\item[$(3\Rightarrow 1)$]
Let $m=[a_1,a_2,d,b,c]$ be a configuration satisfying 3.
$c''b''=0$ implies either $c''=0$ or $b''=0$. It follows that the pair
$({\rm Span}\{(0,1)\}, {\rm Span}\{(0,1)\})$ is a special pair 
hence the configuration is degenerate. 
Now, from proposition \ref{prop1nnX1} it follows that $m$ is equivalent 
to the completely reducible configuration
\[
m'\oplus m''=(a_1',a_2',0,b',c')\oplus(z_1,z_2,0,0)
\]
Notice that $(a_1',a_2',0,b',c')\in S_0\Mm_1(X_{1L})$.
Then, from proposition \ref{prop50}, there is
$\tilde m\in C$ such that $\pi_{R*}\tilde m^{\vee\vee}=m'$.
Then, from the characterization of points in the completion it follows
that $\pi_{R*}\tilde m=m$. 
\end{enumerate}
\end{proof}

The homotopy equivalences (4) and (5) are a direct consequence
of the proposition

\begin{prop}
Let
\begin{align*}
&N_z=\left\{\,(a_{1z},a_{2z},b_z,c_z)\in\Mm_1(X_0)\,|\,|a_{1z}-z|<\delta\,
\right\}\\
&N'=\left\{\,(a_1',a_2',d',b',c')\in\Mm_1(X_1)\,|\,|d'a_1'|<\delta\,
\right\}
\end{align*}
Let $N_0\subset\Mm_2(X_0)$ be the subset of points $(a_1,a_2,b,c)$
with the eigenvalues of $a_1$ lying in $\delta$ neighborhoods of
$x_L$ and $x_R$. Consider the maps
$\boxplus_0:N_{x_{1L}}\times N_{x_{1R}}\to N_0$
defined by
\[
[a_{1L},a_{2L},b_L,c_L]\,\boxplus_0\,[a_{1R},a_{2R},b_R,c_R]
=[a_1,a_2,b,c]
\]
\[
a_1=\begin{bmatrix}a_{1L}&0\\0&a_{1R}\end{bmatrix}\ ,\ 
a_2=\begin{bmatrix}
a_{2L}&\frac{b_Lc_R}{a_{1R}-a_{1L}}\\\frac{b_Rc_L}{a_{1L}-a_{1R}}&a_{2R}
\end{bmatrix}\ ,\ 
b=\begin{bmatrix}b_L\\b_R\end{bmatrix}\ ,\ 
c=\begin{bmatrix}c_L&c_R\end{bmatrix}
\]
and $\boxplus_L: N'\times  N_{z_1}\to  N_L$ defined by
\[
[a_1',a_2',d',b',c']\,\boxplus_L\,[a_1'',a_2'',b'',c'']=
[a_1,a_2,d,b,c]
\]
\begin{gather}\label{eq5R}
a_1=\begin{bmatrix}
a_1'&0\\0&a_1''\end{bmatrix}\ ,\ 
a_2=\begin{bmatrix}
a_2'&\frac{b'c''}{a_1''-d'a_1'}\\\frac{b''c'}{d'a_1'-a_1''}&a_2''
\end{bmatrix}\\
\notag
d=\begin{bmatrix}
d'&0\\0&1\end{bmatrix}\ ,\ 
b=\begin{bmatrix}b'\\b''\end{bmatrix}\ ,\ 
c=\begin{bmatrix}c'&c''\end{bmatrix}
\end{gather}
Then
\begin{enumerate}
\item The maps $\boxplus_0,\boxplus_L$ are homeomorphisms;
\item The inclusions $N_z\to\Mm_1(X_0)$, $N'\to\Mm_1(X_1)$
are homotopy equivalences;
\item $\pi_R^*N_L\cap\pi_L^*N_R=\pi_\emptyset^*N_0$.
\end{enumerate}
\end{prop}
\begin{proof}
Statement (2) is clear from the definition. To prove statement (3)
we observe that 
\[
\pi_R^*N_L\cap\pi_L^*N_R=
\pi_R^*N_L\cap\pi_\emptyset^*\Mm_2(X_0)\cong
N_L\cap\pi_L^*\Mm_2(X_0)
\]
The result now follows easily from proposition \ref{proptau1}.
We turn to the proof of statement (1).
It is an easy consequence of proposition \ref{prop1nnX1} that 
$\boxplus_0$ and $\boxplus_L$ preserve the nondegeneracy of the 
configurations so the maps are well defined.

Now we look at $\boxplus_L$. For $\delta$ small enough the eigenvalues
of $a_1d$ are distinct. Hence
we can choose, up to the action of $(\Cc^*)^{\times4}$,
eigenvector basis $\{v_0,v_1\}\subset V$ of $a_1d$ and
$\{w_0,w_1\}\subset W$ of $da_1$,
where $v_0,w_0$ correspond to the eigenvalues near $0$. 
Normalize $v_1,w_1$ so that $dv_1=w_1$.
Then the action of $(\Cc^*)^{\times4}$ is
reduced to an action of $(\Cc^*)^{\times3}$. We can thus write 
(see also equation (\ref{eq5nboxplusL}))
\begin{gather}
a_1=\begin{bmatrix}
a_1'&0\\0&a_1''\end{bmatrix}\ ,\ 
a_2=\begin{bmatrix}
a_2'&\frac{b'c''}{a_1''-d'a_1'}\\\frac{b''c'}{d'a_1'-a_1''}&a_2''
\end{bmatrix}\\ \notag
d=\begin{bmatrix}
d'&0\\0&1\end{bmatrix}\ ,\ 
b=\begin{bmatrix}b'\\b''\end{bmatrix}\ ,\ 
c=\begin{bmatrix}c'&c''\end{bmatrix}
\end{gather}
The group $(\Cc^*)^{\times3}$ acts transitively
on equivalence classes of such configurations
written in the above canonical form. 
This shows the existence of an inverse, 
hence $\boxplus_L$ is a homeomorphism. The proof for $\boxplus_0$ is similar.
\end{proof}

The maps $\boxplus_0,\boxplus_L$ extend to the closure 
$\bar N',\bar N_z$ of $N',N_z$.
The following proposition is a direct consequence of proposition
\ref{prop1nnX1}:

\begin{prop}\label{prop53}\label{prop5ndegenerateL}$\,$
\begin{itemize}
\item Let $m_L=[a_{1L},a_{2L},b_L,c_L]\in\bar N_{x_{1L}}$,
$m_R=[a_{1R},a_{2R},b_R,c_R]\in\bar N_{x_{1R}}$.
Then the following are equivalent:
\begin{enumerate}
\item $m_L\boxplus_0m_R$ is degenerate;
\item Either $m_L$ or $m_R$ is degenerate.
\item At least one of the 4 vectors
$b_L,b_R,c_L,c_R$ is zero.
\end{enumerate}
\item Let $m'=[a_1',a_2',d',b',c']\in \bar N'$, 
$m''=[a_1'',a_2'',b'',c'']\in\bar N_{z_1}$.
The following are equivalent:
\begin{enumerate}
\item $m'\boxplus_Lm''$ is degenerate;
\item Either $m'$ or $m''$ is degenerate;
\item One of the 4 vectors $b',b'',c',c''$ is zero.
\end{enumerate}
\end{itemize}
\end{prop}

We are ready to prove

\begin{prop}
$N_2$ is an open neighborhood of $C$.
\end{prop}
\begin{proof}
From lemma \ref{lemm511} it follows immediately
that $\pi_{R*}C\subset\bar N_L$.
\item
Suppose there is a sequence $y_n\in\Mm_2(X_{1L})$
such that $y_n\to y\in\pi_{R*}C$. Write $y_n=[a_{1n},a_{2n},d_n,b_n,c_n]$.
Then, by property 2 in lemma \ref{lemm511}
the eigenvalues of $d_na_{in}$ converge to $0,z_i$. Hence, for $n$
large enough $y_n\in N_L$. Hence $N_L\cup\pi_{R*}C$ is an open
neighborhood of $\pi_{R*}C$.
\item
Suppose there is a sequence $x_n\to x\in C$ such that
$x_n\notin N_2$. Hence $x_n\notin C$ so,
by passing to a subsequence we may assume without loss of
generality that $x_n\in\pi_R^*\Mm_2(X_{1L})$.
Let $y_n=\pi_{R*}x_n\in\Mm_2(X_{1L})$ and
write $y_n=[a_{1n},a_{2n},d_n,b_n,c_n]$.
Then $y_n\to y=\pi_{R*}x$
by continuity of $\pi_{R*}$,
and $y_n\notin N_L$.
But by property 2 in lemma \ref{lemm511}
the eigenvalues of $d_na_{in}$ converge to $0,z_i$ which implies,
for $n$ large enough, that $y_n\in N_L$.
\end{proof}

Finaly we prove the homotopy equivalence (6):

\begin{prop}\label{prop(6)}
The inclusion $C\to N_2$ is a strong deformation retract.
\end{prop}
\begin{proof}
We will construct a homotopy $H_2:N_2\times[0,1]\to N_2$ betwen
the identity and a retraction $N_2\to C$.
Let $H_{x_1,x_2}:\bar N_z\times[0,1]\to \bar N_z$ be defined by
\[
H_{x_1,x_2}(a_1,a_2,b,c,t)=
\left(\,t^2a_{1}+(1-t^2)x_{1},t^2a_{2}+(1-t^2)x_{2},tb,tc\,\right)
\]
and let $H_1:\bar N'\times[0,1]\to\bar N'$ be defined by
\[
H_1(a_1',a_2',d',b',c',t)=(a_1',a_2',t^2d',b',c')
\]
Then we defined $H_L:\bar N_L\times[0,1]\to\bar N_L$ by
\[
H_L(m'\boxplus_Lm'',t)\defeq H_1(m',t)\boxplus_LH_{z_1,z_2}(m'',t)
\]
We define $H_2$ as the unique solution of
the system of equations
\begin{gather}\label{definingH2}
\pi_{R*}H_2(x,t)=H_L(\pi_{R*}x,t)\\\notag
\pi_{L*}H_2(x,t)=H_R(\pi_{L*}x,t)
\end{gather}
We have to show existence and uniqueness of solution.
Then we will show that $H_2$
defines a homotopy between the identity on $N_2$ and a retraction
$N_2\to C$.

We define the auxiliary map
$H_0:\bar N_0\times[0,1]\to\bar N_0$ by
\[
H_0(m_L\boxplus_0m_R,t)\defeq
H_{x_{1L},x_{2L}}(m_L,t) 
\boxplus_0 
H_{x_{1R},x_{2R}}(m_R,t)
\]
To prove existence and uniqueness of solution of the
system (\ref{definingH2}) we consider two cases:
\begin{enumerate}
\item Assume that either $t=0$ or $x\in C$. Then we claim that
$H_L(\pi_{R*}x,t)\in\pi_{R*}C$, $H_R(\pi_{L*}x,t)\in\pi_{L*}C$.
If $t=0$ this follows from directly from lemma \ref{lemm511}.
If $x\in C$ then, from lemma \ref{lemm511}
we can write
\[
\pi_{R*}x=x'\boxplus_L x''=(a_1',a_2',0,b',c')\boxplus_L(a_1'',a_2'',b'',c'')
\]
with $c''b''=0$. It then follows from the definition of $H_L$
that $H_L(\pi_{R*}x,t)=\pi_{R*}x$ for all $t$. 
In the same way we see that
$H_R(\pi_{L*}x,t)=\pi_{L*}x$. This proves the claim.
Then, existence and uniqueness follows from proposition \ref{prop50}.
\item Assume $t\neq1$ and $x\notin C$. Then we may assume
$\pi_{R*}x\in N_L$. Then, since $H_L(\pi_{R*}x,t)\in N_L$,
we get from (\ref{definingH2})
\[
\pi_{R*}H_2(x,t)=H_L(\pi_{R*}x,t)\Rightarrow H_2(x,t)=\pi_R^*H_L(\pi_{R*}x,t)
\]
This proves uniqueness. To prove existence we need to show that
\[
\pi_{L*}H_2(x,t)=\pi_{L*}\pi_R^*H_L(\pi_{R*}x,t)=H_R(\pi_{L*}x,t)
\]
It is enough to show this for the case where
$x=\pi_L^*\pi_R^*y$ for some $y\in N_0$ since the set of points of this
form is dense and $H_L,H_R,\pi_{L*},\pi_R^*,\pi_{R*}$ are continuous.

It is an easy computation to show that
$H_L(\pi_L^*y,t)=\pi_L^*H_0(y,t)$, $H_R(\pi_R^*y,t)=\pi_R^*H_0(y,t)$. 
It follows that
\begin{multline*}
\ \ \ \ \pi_{L*}\pi_R^*H_L(\pi_{R*}x,t)=\pi_{L*}\pi_R^*\pi_L^*H_0(y,t)=\\
=\pi_R^*H_0(y,t)=H_R(\pi_{L*}x,t)
\end{multline*}
\end{enumerate}
Now we need to show $H_2$ is the desired homotopy.
 Direct inspection shows $H_2(x,1)=x$.
We saw in (1) above that, for $x\in C$, $H_2(x,t)=x$ and
$H_2(x,0)\in C$.
The continuity of $H_2$ follows from the continuity of
$\pi_{L*},\pi_{R*},H_L,H_R$.
\end{proof}

\section{The differential $d_1$}

The objective of this section is to
obtain the homotopy type of the inclusion maps
\[
\xymatrix{
A_L&A_0\ar[l]\ar[r]&A_R\\
N_L\ar[u]\ar[rd]&N_0\ar[u]\ar[l]\ar[r]&N_R\ar[u]\ar[ld]\\
&N_2
}\]
where $A_0=A_L\cap A_R=\pi_\emptyset^*\Mm_2(X_0)$.
Since these spaces are classifying spaces it is enough to study
the pullback under these maps of the tautological bundles.
Together with the open cover of the previous section this will
give a description of the homotopy type of $\Mm_2(X_2)$. It will
also allow us to
compute the $d_1$ differential in the spectral sequence introduced
in corolary \ref{corospseq2}.

\begin{lema}
Let
\begin{gather*}
\tilde F_0(k,r)=\frac{
\left\{(u,v):u,v:\Cc^k\to \Cc^r\,,\,u,v\mbox{ are injective}\right\}}
{(u,v)\sim\left(u(\bar g^t)^{-1},vg\right),\,g\in Gl(k,\Cc)}\\
\tilde F_1(k,r)=\frac{
\left\{(u,v):u,v:\Cc^k\to \Cc^r\,,\,u,v\mbox{ are injective}\right\}}
{(u,v)\sim\left(u(\bar g_u^t)^{-1},vg_v\right),\,g_u,g_v\in Gl(k,\Cc)}
\end{gather*}
and define $F_0(k,r)\subset \tilde F_0(k,r)$ and
$F_1(k,r)\subset \tilde F_1(k,r)$ by
\begin{gather*}
F_0(k,r)=\left\{\left[\bar b^t,c\right]\in\tilde
F_0(k,r)\,,\,bc=0\right\}\\
F_1(k,r)=\left\{\left[\bar b^t,c\right]\in\tilde
F_1(k,r)\,,\,bc=0\right\}
\end{gather*}
Let
$j_0:F_0(k,r)\to\Mm_k^r(X_0)$,
$j_1:F_1(k,r)\to\Mm_k^r(X_1)$ be the inclusion maps given by
$\left[\bar b^t,c\right]\mapsto[b,c]$. Then
we have the homotopy commutative diagram
\begin{equation}\label{diagipj}
\xymatrix{
&\Mm_k^r(X_0)\ar[r]^-{\pi^*}&\Mm_k^r(X_1)\\
Gr(k,\Cc^r)&
F_0(k,r)\ar[l]^-{p_0}\ar[u]^-{j_0}\ar[r]^-{pr}\ar[d]^-{\imath_0}&
F_1(k,r)\ar[u]^-{j_1}\ar[d]^-{\imath_1}\\
&\tilde F_0(k,r)\ar[lu]^-{\tilde p_0}\ar[r]^-{\widetilde{pr}}&
\tilde F_1(k,r)
}\end{equation}
where $p_0$ is the projection $[b,c]\mapsto[c]$.
Moreover, in the rank stable limit,
the maps $\imath_0,\imath_1,p_0,\tilde p_0,j_0,j_1$
are homotopy equivalences. 
\end{lema}
\begin{proof}
We divide the proof into three steps:
\begin{enumerate}
\item[Step 1.]
$p_0,\tilde p_0$ are fibrations with fibers $M(k,r-k)$
and $M(k,r)$ respectively where
$M(k,r)=\frac{U(r)}{U(k)}$ is the space of injective maps
from $\Cc^k$ to $\Cc^{r}$ which is contractible in the stable range.
That proves $p_0,\tilde p_0$ are homotopy equivalences.
It imediatelly follows that $\imath_0$ is a homotopy equivalence.
\item[Step 2.] Now we look at $\imath_1$. 
Consider the projection $p_1:F_1(k,r)\to Gr(k,\Cc^r)$
given by $[b,c]\mapsto [c]$. When $r\to\infty$,
the spectral sequence associated with the fibration
\[
\xymatrix{
Gr(k,\Cc^{r-k})\ar[r]^-{i}&F_1(k,r)\ar[d]^-{p_1}\\
&Gr(k,\Cc^r)
}\]
collapses since all homology is in even dimensions. 
It easilly follows that $\imath_1$ is an isomorphism in all
homology groups, hence an homotopy equivalence.
\item[Step 3.]
Finally we need to prove the statements about $j_0,j_1$. Let
$C_k^r(X_0)$, $C_k^r(X_1)$ be the spaces of configurations
corresponding to the monads for $X_0$ and $X_1$.
Let $C_0^F(k,r)\subset C_k^r(X_0)$ and
$C_1^F(k,r)\subset C_k^r(X_1)$ be the subsets of
configurations of the form $(0,0,b,c)$ and
$(0,0,0,b,c)$ respectively. Then we have the maps
between fibrations
\[
\xymatrix{
Gl(k)\ar[r]\ar@{=}[d]&C_0^F(k,r)\ar[d]\ar[r]&F_0(k,r)\ar[d]^{j_0}\\
Gl(k)\ar[r]&C_k^r(X_0)\ar[r]&\Mm_k^r(X_0)
}\]
and a similar diagram for $X_1$. In the rank stable limit
the spaces $C_0^F(k,r)$ and $C_k^r(X_0)$ are contractible
(see \cite{San95}, \cite{BrSa97})
so, by the five lemma $j_0$ is an isomorphism in $\pi_*$ hence
an homotopy equivalence. The proof for $j_1$ is the same.
\end{enumerate}
\end{proof}

We are ready to state and prove the main theorem of this section:

\begin{teor}\label{teor6main}
Consider the compositions
\begin{align}
&\xymatrix{
A_0\ar[r]^-{\pi_{\emptyset*}}&\Mm_2^\infty(X_0)\ar[r]^-{j_0^{-1}}&
F_0(2,\infty)\ar[r]^-{p_0}&
Gr(2,\Cc^\infty)}\label{align5nA0}\\
&\xymatrix{
A_L\ar[r]^-{\pi_{R*}}&\Mm_2^\infty(X_{1L})\ar[r]^-{j_1^{-1}}&
F_1(2,\infty)\ar[r]^-{\imath_1}&
\tilde F_1(2,\infty)}\label{align5nAL}\\
&\xymatrix{
N_0\ar[r]^-{\boxplus_0^{-1}}&
N_{x_{1L}}\times N_{x_{1R}}\ar[r]^-{p_L}&
N_{x_{1L}}\ar[r]&
\Mm_1(X_0)\ar[r]^-{p_0j_0^{-1}}&
Gr(1,\Cc^\infty)}\label{align5nN0L}\\
&\xymatrix{
N_0\ar[r]^-{\boxplus_0^{-1}}&
N_{x_{1L}}\times N_{x_{1R}}\ar[r]^-{p_R}&
N_{x_{1R}}\ar[r]&
\Mm_1(X_0)\ar[r]^-{p_0j_0^{-1}}&
Gr(1,\Cc^\infty)}\label{align5nN0R}\\
&\xymatrix{
N_L\ar[r]^-{\boxplus_L^{-1}}&
N'\times N_{z_1}\ar[r]^-{p''}&
N_{z_1}\ar[r]&
\Mm_1(X_0)\ar[r]^-{p_0j_0^{-1}}&
Gr(1,\Cc^\infty)}\label{align5nNL0}\\
&\xymatrix{
N_L\ar[r]^-{\boxplus_L^{-1}}&
N'\times N_{z_1}\ar[r]^-{p'}&
N'\ar[r]&
\Mm_1(X_1)\ar[r]^-{\imath_1j_1^{-1}}&
\tilde F_1(1,\infty)}\label{align5nNL1}\\
&\xymatrix{
N_2\ar[r]^-{\simeq}&
C\ar[r]^-{\pi_{L*}^{\vee\vee}}&
S_0\Mm_1^\infty(X_{1L})\ar@{ >->}[r]&
\Mm_1(X_{1L})\ar[r]^-{\imath_1j_1^{-1}}&
\tilde F_1(1,\infty)}\label{align5nCL}\\
&\xymatrix{
N_2\ar[r]^-{\simeq}&
C\ar[r]^-{\pi_{R*}^{\vee\vee}}&
S_0\Mm_1^\infty(X_{1R})\ar@{ >->}[r]&
\Mm_1(X_{1R})\ar[r]^-{\imath_1j_1^{-1}}&
\tilde F_1(1,\infty)}\label{align5nCR}
\end{align}
Let $E,L$ be the tautological bundles over $Gr(2,\infty)$
and $Gr(1,\infty)$ respectivelly. Then we define
the following bundles:
\begin{itemize}
\item $E_0\to A_0$ is the pullback of $E$
under the composition \ref{align5nA0}.
\item $L_{0L,0}\to N_0$ is the pullback of $L$ under \ref{align5nN0L}
\item $L_{0R,0}\to N_0$ is the pullback of $L$ under \ref{align5nN0R}
\item $L_{0R,L}\to N_L$ is the pullback of $L$ under \ref{align5nNL0}
\end{itemize}
Now let $\tilde E_u,\tilde E_v\to\tilde F_1(2,r)$ be the tautological
bundles corresponding to $u,v$ and let
$\tilde L_u,\tilde L_v$ be the tautological 
line bundles over $\tilde F_1(1,\infty)$. We define
\begin{itemize}
\item $E_{bL},E_{cL}\to A_L$ are the pullback of 
$\tilde E_u,\tilde E_v$ under \ref{align5nAL}.
\item $L_{bL,L},L_{cL,L}\to N_L$ are the pullback of 
$\tilde L_u,\tilde L_v$ under \ref{align5nNL1}
\item $L_{bL,2},L_{cL,2}\to N_2$ are the pullback of
$\tilde L_u,\tilde L_v$ under \ref{align5nCL}
\item $L_{bR,2},L_{cR,2}\to N_2$ are the pullback of
$\tilde L_u,\tilde L_v$ under \ref{align5nCR}
\end{itemize}
Then we have
\begin{enumerate}
\item[\bf (1)]
$\displaystyle E_{bL}|_{A_0}=E_0\ ,\ E_{cL}|_{A_0}=E_0$
\vspace*{.2em}
\item[\bf (2)]
$\displaystyle
L_{bL,L}|_{N_0}\cong L_{0L,0}\ ,\ L_{cL,L}|_{N_0}\cong L_{0L,0}\ ,\ 
L_{0R,L}|_{N_0}\cong L_{0R,0}$
\vspace*{.2em}
\item[\bf (3)]
$\displaystyle
E_{bL}|_{N_L}\cong L_{bL}\oplus L_{0R}\ ,\ 
E_{cL}|_{N_L}\cong L_{cL}\oplus L_{0R}$
\vspace*{.2em}
\item[\bf (4)]
$\displaystyle E_0|_{N_0}\cong L_{0L,0}\oplus L_{0R,0}$.
\vspace*{.2em}
\item[\bf (5)]
$\displaystyle L_{bR,2}|_{N_L}\cong L_{cR,2}|_{N_L}\cong L_{0R,L}$
\vspace*{.2em}
\item[\bf (6)]
$\displaystyle
L_{bL,2}|_{N_L}\cong L_{bL,L}\ ,\ L_{cL,2}|_{N_L}\cong L_{cL,L}$
\end{enumerate}
Similar statements hold for the spaces $A_R,N_R$ and the maps
$N_R\to A_R$, $N_R\to N_2$ and $N_0\to N_R$.
\end{teor}
\begin{proof}$\,$
\begin{enumerate}
\item[\bf (1)] First we show that $E_{bL}|_{A_0}\cong E_{cL}|_{A_0}\cong E_0$.
Consider diagram (\ref{diagipj}). We will start by defining a homotopy
inverse $q:Gr(k,\Cc^r)\to\tilde F_0(k,r)$
to the map $\tilde p_0:\tilde F_0\to Gr(k,\Cc^\infty)$ as follows:
choose a map $c:\Cc^k\to\Cc^r$ representing an element
$[c]\in Gr(k,\Cc^r)$. Choose $h\in Gl(k,\Cc)$
such that $ch$ is orthogonal.
Then define $q([c])=[ch,ch]$.
This map is well defined and independent of the choice of $h$.
Also $p_0q=\Id$ hence $p_0=q^{-1}$.

Now observe that the composition
\[
\tilde{pr}\circ q:Gr(k,\Cc^r)\to\tilde F_1(k,r)=Gr(k,\Cc^r)\times
Gr(k,\Cc^r)
\]
is the diagonal map. It follows that, if
$E$ is the tautological bundle over $Gr(k,\Cc^\infty)$, then
\[
q^*\tilde{pr}^*\tilde E_u\cong q^*\tilde{pr}^*\tilde E_v\cong E
\]
To show that $E_{bL}|_{A_0}\cong E_{cL}|_{A_0}\cong E_0$
it suffices to show that 
$pr^*\imath_1^*\tilde E_u\cong pr^*\imath_1^*\tilde E_v\cong p_0^*E$.
We have
\[
pr^*\imath_1^*\tilde E_u=\imath_0^*\widetilde{pr}^*\tilde E_u\cong
p_0^*q^*\widetilde{pr}^*\tilde E_u=p_0^*E
\]
and a similar statement is true for $\tilde E_v$. This concludes
the proof.
\item[\bf (2)] We want to show that
\[
L_{bL,L}|_{N_0}\cong L_{0L,0}\ ,\ L_{cL,L}|_{N_0}\cong L_{0L,0}\ ,\ 
L_{0R,L}|_{N_0}\cong L_{0R,0}
\]
We have the commutative diagram (see proposition \ref{proptau2})
\[
\xymatrix{
N_0\ar[d]^{\pi^*}&
N_{x_{1L}}\times N_{x_{1R}}
\ar[l]_-{\boxplus_0}\ar[r]^-{p_R}\ar[d]^{\pi^*\times\tau}&
N_{x_{1R}}\ar[r]\ar[d]^{\tau}&
\Mm_1(X_0)\ar[d]^{\tau}&
F_0\ar[l]_-{j_0}\ar[r]^-{p_0}\ar[d]&
Gr\ar[d]\\
N_L&
N'\times N_z\ar[l]_-{\boxplus_L}\ar[r]^-{p''}&
N_z\ar[r]&
\Mm_1(X_0)&
F_0\ar[l]_-{j_0}\ar[r]^-{p_0}&
Gr
}\]
from which it follows that $L_{0R,L}|_{N_0}\cong L_{0R,0}$.
We also have the commutative diagram
\[
\xymatrix{
&&&&&Gr\\
N_0\ar[d]^{\pi^*}&
N_{x_{1L}}\times N_{x_{1R}}
\ar[l]_-{\boxplus_0}\ar[r]^-{p_L}\ar[d]^{\pi^*\times\tau}&
N_{x_{1L}}\ar[r]\ar[d]^{\tau}&
\Mm_1(X_0)\ar[d]^{\tau}&
F_0\ar[l]_-{j_0}\ar[r]^-{\imath_0}\ar[d]^{pr}\ar[ru]^-{p_0}&
\tilde F_0\ar[u]_{\tilde p_0}\ar[d]^{\tilde{pr}}\\
N_L&
N'\times N_z\ar[l]_-{\boxplus_L}\ar[r]^-{p'}&
N'\ar[r]&
\Mm_1(X_1)&
F_1\ar[l]_-{j_1}\ar[r]^-{\imath_1}&
\tilde F_1
}\]
from which it follows, 
as in step (1) above, that
$L_{bL,L}|_{N_0}\cong L_{cL,L}|_{N_0}\cong L_{0L,0}$.
\item[\bf (3)] We want to show that
\[
E_{bL}|_{N_L}\cong L_{bL,L}\oplus L_{0R,L}\ ,\ 
E_{cL}|_{N_L}\cong L_{cL,L}\oplus L_{0R,L}
\]
Consider the following diagram:
\begin{equation}\label{diag5nNL->AL}
\xymatrix{
A_L&&
F_1\ar[ll]_-{j_1}^-{\simeq}\ar[rr]^-{\imath_1}_-{\simeq}&&
\tilde F_1\\
N_L\ar[u]&
N'\times N_{z_1}\ar[l]_-{\boxplus_L}^-{\cong}\ar[r]_-{\simeq}&
\Mm_1(X_1)\times\Mm_1(X_0)&
F_1\times F_0\ar[l]_-{j_1\times j_0}^-{\simeq}
\ar[r]^-{\imath_1\times \imath_0}_{\simeq}&
\tilde F_1\times \tilde F_0\ar@{-->}[u]^{\tilde w}
}\end{equation}
Since 
$\tilde{pr}^*\tilde L_{u}\cong\tilde{pr}^*\tilde L_{v}\cong \tilde p_0^*L$,
the proof will be complete if we show there is a map 
$\tilde w:\tilde F_1(1,\infty)\times \tilde F_0(1,\Cc^\infty)\to
\tilde F_1(2,\infty)$ making the diagram homotopy commutative, such that
\begin{equation}\label{eq5nnw}
\tilde w^*\tilde E_u=\tilde L_{u}\oplus \tilde{pr}^*\tilde L_{u}\ ,\ 
\tilde w^*\tilde E_v=\tilde L_{v}\oplus \tilde{pr}^*\tilde L_{v}
\end{equation}
We begin by building $\tilde w$. 
Define maps $s_L,s_R:Gr(1,\Cc^\infty)\to Gr(1,\Cc^\infty)$
as follows:
let $v:\Cc\to\Cc^\infty$ and write $v=(v^1,v^2,\ldots)$. Then
\begin{align*}
&s_L([v])\defeq[(v^1,0,v^2,0,\ldots)]\\
&s_R([v])\defeq[(0,v^1,0,v^2,\ldots)]
\end{align*}
We observe that $s_L,s_R$ are homotopic to the identity. It follows that,
if we define
\[
\tilde w:([b_L,c_L],[b_R,c_R])\mapsto 
[s_L(b_L)\oplus s_R(b_R),s_L(c_L)\oplus s_R(c_R)]
\]
then
\[
(\tilde w)^*\tilde E_u=
\tilde L_{u}\oplus\tilde{pr}^*\tilde L_{u}\ ,\ 
(\tilde w)^*\tilde E_v=
\tilde L_{v}\oplus \tilde{pr}^*\tilde L_{v}
\]
It remain to show diagram \ref{diag5nNL->AL} is commutative.
Let $j_z:F_0(1,\infty)\to N_z$ be defined by $j_z:[b,c]\mapsto[z,0,b,c]$.
Then the diagram
\[
\xymatrix{
N'\times N_{z_1}\ar[rd]&&
F_1\times F_0\ar[ll]_-{j_1\times j_{z_1}}\ar[ld]\\
&\Mm_1(X_1)\times\Mm_1(X_0)
}\]
is homotopy commutative. 
We are left with the diagram
\[
\xymatrix{
A_L&&
F_1\ar[ll]_-{j_1}^-{\simeq}\ar[rr]^-{\imath_1}_-{\simeq}&&
\tilde F_1\\
N_L\ar[u]&
N'\times N_{z_1}\ar[l]_-{\boxplus_L}^-{\cong}&&
F_1\times F_0\ar[ll]_-{j_1\times j_{z_1}}^-{\simeq}
\ar[r]^-{\imath_1\times \imath_0}_{\simeq}&
\tilde F_1\times \tilde F_0\ar[u]^{\tilde w}
}\]
Now define the map $w:F_1(1,\infty)\times F_0(1,\infty)\to F_1(2,\infty)$ by
\[
w:([b_L,c_L],[b_R,c_R])\mapsto 
[s_L(b_L)\oplus s_R(b_R),s_L(c_L)\oplus s_R(c_R)]
\]
Clearly we have the commutative diagram
\[
\xymatrix{
F_1(2,\infty)\ar[rr]^-{\imath_1}&&
\tilde F_1(2,\infty)\\
F_1(1,\infty)\times F_0(1,\infty)\ar[rr]^-{\imath_1\times\imath_0}
\ar[u]^{w}&&
\tilde F_1(1,\infty)\times\tilde F_0(1,\infty)\ar[u]^{\tilde w}
}\]
We are thus left with the diagram
\[
\xymatrix{
A_L&&&
F_1\ar[lll]_-{j_1}^-{\simeq}\\
N_L\ar[u]&
N'\times N_{z_1}\ar[l]_-{\boxplus_L}^-{\cong}&&
F_1\times F_0\ar[ll]_-{j_1\times j_{z_1}}^-{\simeq}\ar[u]^w
}\]
Next we introduce maps
\begin{align*}
&S_L([b,c])=[s_L(b),s_L(c)]\\
&S_R([b,c])=[s_R(b),s_R(c)]
\end{align*}
These maps are homotopic to the identity hence we only have to show the
diagram
\[
\xymatrix{
A_L&&&&&
F_1\ar[lllll]_-{j_1}\\
N_L\ar[u]&
N'\times N_{z_1}\ar[l]_-{\boxplus_L}&&
F_1\times F_0\ar[ll]_-{j_1\times j_{z_1}}&&
F_1\times F_0\ar[ll]_-{S_L\times S_R}\ar[u]^w
}\]
is homotopy commutative. 
This is an easy direct verification.
\item[\bf (4)] We want to show that
$E_0|_{N_0}\cong L_{0L}\oplus L_{0R}$.
Consider the diagram
\[
\xymatrix{
N_0\ar[r]^{i_1}\ar[d]^{i_2}&A_0\ar[d]^{i_3}\\N_L\ar[r]^{i_4}&A_L
}\]
Then $E_0=i_3^*E_{bL}$ so
\[
E_0|_{N_0}=i_1^*E_0=i_1^*i_3^*E_{bL}=i_2^*i_4^*E_{bL}=L_{0L}\oplus L_{0R}
\]
\item[\bf (5)] We want to show that
$L_{bR}|_{N_L}\cong L_{cR}|_{N_L}\cong L_{0R}$.
The result will follow if we show that the following diagram is homotopy
commutative:
\begin{equation}\label{diag5nNL->N2}
\xymatrix{
N_L\ar[d]&
N'\times N_{z_1}\ar[l]_-{\boxplus_L}^-{\cong}\ar[r]&
N_{z_1}\ar[r]&
\Mm_1(X_0)\ar[d]^{\pi^*}&
F_0\ar[l]\ar[r]&\tilde F_0\ar[d]^{\tilde{pr}}\\
N_2&
C\ar[l]_-{\simeq}\ar[r]^-{\pi_{L*}^{\vee\vee}}&
S_0\Mm_1(X_1)\ar[r]^-{\simeq}&
\Mm_1(X_1)&
F_1\ar[l]\ar[r]&\tilde F_1
}\end{equation}
Let
\begin{align*}
&S_{1L}N_2=\left\{(\Ee,\phi)\in N_2\,|\,
c_2\left((\pi_{L*}\Ee)^{\vee\vee}\right)=1\right\}\\
&S_1N_L=\left\{(\Ee,\phi)\in N_L\,|\,
c_2\left((\pi_{L*}\Ee)^{\vee\vee}\right)=1\right\}\\
&S_0N'=\left\{(\Ee,\phi)\in N'\,|\,
c_2\left((\pi_{L*}\Ee)^{\vee\vee}\right)=0\right\}
\end{align*}
Then the commutativity of diagram (\ref{diag5nNL->N2})
follows from the commutativity of
\begin{equation}\label{diag5nnNL->N2big}
\xymatrix{
N_2&&N_L\ar[ll]^-{\pi_R^*}&&N'\times N_{z_1}\ar[ll]_-{\boxplus_L}\ar[dd]\\
C\ar@{ >->}[r]\ar[d]_-{\pi_{L*}^{\vee\vee}}\ar@{ >->}[u]&
S_{1L}N_2\ar@{ >->}[lu]\ar[ld]^{\pi_{L*}^{\vee\vee}}&
S_1N_L\ar[l]^-{\pi_R^*}\ar@{ >->}[u]\ar[d]^{\pi_{L*}^{\vee\vee}}&
S_0N'\times N_{z_1}\ar@{ >->}[ru]\ar[l]_-{\boxplus_L}\ar@{->>}[rd]\\
\Mm_1(X_{1R})&&\Mm_1(X_0)\ar[ll]_-{\pi_R^*}&&N_{z_1}\ar@{ >->}[ll]
}\end{equation}
We need to check the image of $\boxplus_L:S_0N'\times N_{z_1}\to N_L$
is contained in $S_1N_L$. Then,
analyzing the commutativity of diagram (\ref{diag5nnNL->N2big})
boils down do analyzing the diagram
\begin{equation}\label{diag5nnNL->N2small}
\xymatrix{
S_1N_L\ar[d]^{\pi_{L*}^{\vee\vee}}&
S_0N'\times N_{z_1}\ar[l]_-{\boxplus_L}\ar[d]\\
\Mm_1(X_0)&N_{z_1}\ar[l]
}\end{equation}
Let $m'\in S_0N'\subset S_0\Mm_1(X_1)$,
$m'=[a_1',a_2',0,b',c']$. 
Let $m''\in N_z$. Then a direct computation shows that
$(\pi_{L*}(m'\boxplus_Lm''))^{\vee\vee}=m''$.
This shows that the image of $S_0N'\times N_{z_1}$ under $\boxplus_L$
is contained in $S_1N_L$ and that diagram (\ref{diag5nnNL->N2small})
is commutative.
\item[\bf (6)] We want to show that $L_{bL,2}|_{N_L}\cong L_{bL,L}$,
$L_{cL,2}|_{N_L}\cong L_{cL,L}$. This will follow from the commutativity
of the diagram
\begin{equation}\label{diag5bnNL->N2}
\xymatrix{
N_L\ar[d]^{\pi_R^*}&
N'\times N_{z_1}\ar[l]_-{\boxplus_L}^-{\cong}\ar[r]&
N'\ar[r]&
\Mm_1(X_1)\ar@{=}[d]\\
N_2&
C\ar[l]_-{\simeq}\ar[r]^-{\pi_{R*}^{\vee\vee}}&
S_0\Mm_1(X_1)\ar[r]^-{\simeq}&
\Mm_1(X_1)
}\end{equation}
We showed in proposition \ref{prop(6)} that
the map $H_2(\cdot,0)$ is the homotopy inverse of the inclusion $C\to N_2$.
Let $(m',m'')\in N'\times N_z$. Then, by definition of $H_2$,
\[
\pi_{R*}H_2(\pi_R^*(m'\boxplus_L m''),0)=H_L(m'\boxplus_L m'',0)=
H_1(m',0)\boxplus H_z(m'',0)
\]
Hence the diagram
\[
\xymatrix{
N_L\ar[d]^{\pi_R^*}&
N'\times N_{z_1}\ar[l]_-{\boxplus_L}^-{\cong}\ar[r]&
N'\ar[d]^{H_1(\cdot,0)}\\
N_2\ar[r]^-{H_2(\cdot,0)}&
C\ar[r]^-{\pi_{R*}^{\vee\vee}}&
S_0\Mm_1(X_1)
}\]
is commutative. From here it follows easily that diagram (\ref{diag5bnNL->N2})
is commutative.
\end{enumerate}
\end{proof}

\section{The cohomology of $\Mm_2(X_q)$}

The objective of this section is to prove theorem \ref{theo02}.
We begin by proving it for the special case $q=2$:

\begin{teor}\label{theo5nnnk=2}
There is an exact sequence
\[
0\to K_C\to H^*(\Mm_2(X_2))\to H^*(A_L)\oplus H^*(A_R)\to H^*(A_0)\to 0
\]
where $K_C={\rm Ker}\left(\,H^*(C)\to H^*(N_L)\oplus H^*(N_R)\,\right)$.
This sequence splits and we get
\[
 H^*(\Mm_2(X_2))\cong K_C\oplus
{\rm Ker}\left(\,H^*(A_L)\oplus H^*(A_R)\to H^*(A_0)\,\right)
\]
\end{teor}
\begin{proof}
Recall corolary \ref{corospseq2}. We will use this spectral sequence
to compute $H_*(\Mm_2(X_2))$.
Clearly the map $d_1:E_{1,n}\to E_{2,n}$ is surjective hence 
$E_2^{2,n}=0$. Also we notice that $E_1^{p,2n+1}=0$ for any $p$.
It follows that the spectral sequence collapses at the term $E_2$.
We get then
\begin{gather}
\label{eq5nnH2n}
H^{2n}(\,\Mm_2(X_2)\,)=E_\infty^{0,2n}=
{\rm Ker\,}\left(d_1:E_1^{0,2n}\to E_1^{1,2n}\right)\\
\label{eq5nnH2n+1}
H^{2n+1}(\,\Mm_2(X_2)\,)=E_\infty^{1,2n}=
\frac{{\rm Ker}\,\left(d_1:E_1^{1,2n}\to E_1^{2,2n}\right)}{
{\rm Im}\,\left(d_1:E_1^{0,2n}\to E_1^{1,2n}\right)}
\end{gather}
When performing calculations we will use the following sign conventions:
\begin{equation}\label{diag5nnsigns}
\xymatrix{
A_L&A_0\ar[l]_+\ar[r]^+&A_R\\
N_L\ar[u]^{-}\ar[rd]_+&
N_0\ar[u]^+\ar[l]_+\ar[r]^+&
N_R\ar[u]_{-}\ar[ld]^{-}\\
&N_2
}\end{equation}
We begin by defining the following generators of the cohomology 
of $E^{0,2n}_1$:
\begin{align*}
&a_{\Delta L}^i=c_i(E_{cL})-c_i(E_{bL})&
&a_{\Delta R}^i=c_i(E_{cR})-c_i(E_{bR})\\
&a_{bL}^i=c_i(E_{bL})&&a_{bR}^i=c_i(E_{bR})\\
&c_{\Delta L}=c_1(L_{cL})-c_1(L_{bL})&
&c_{\Delta R}=c_1(L_{cR})-c_1(L_{bR})\\
&c_{bL}=c_1(L_{bL})&&c_{bR}=c_1(L_{bR})
\end{align*}
We do the same for $E^{1,2n}_1$:
\begin{align*}
&n_{\Delta L}=c_1(L_{cL})-c_1(L_{bL})&
&n_{\Delta R}=c_1(L_{cR})-c_1(L_{bR})\\
&n_{bL}=c_1(L_{bL})&&n_{bR}=c_1(L_{bR})\\
&n_{0R}=c_1(L_{0R})&&n_{0L}=c_1(L_{0L})\\
&a^i=c_i(E_0)
\end{align*}
and for $E_1^{2,n}$:
\begin{align*}
&n_{0R}=c_1(L_{0R})&&n_{0L}=c_1(L_{0L})
\end{align*}
Then, from theorem \ref{teor6main} it follows that
the map $d_1:E_1^{0,2n}\to E_1^{1,2n}$ may be represented by the
following diagram, where the entries correspond to those in diagram
(\ref{diag5nnsigns}):
\[
\hspace*{-3em}
\xymatrix@C=-4.5em{
\left(a^1_{\Delta L},a^1_{bL},a^2_{\Delta L},a^2_{bL}\right)
\ar@{|->}[r]\ar@{ |->}[d]&
\left(0,a^1,0,a^2\right)&
\left(a^1_{\Delta R},a^1_{bR},a^2_{\Delta R},a^2_{bR}\right)
\ar@{|->}[l]\ar@{ |->}[d]
\\
{\begin{array}{c}
\left(-n_{\Delta L},-n_{bL}-n_{0R},-n_{\Delta L}n_{0R},-n_{bL}n_{0R}\right)\\
\left(n_{\Delta L},n_{bL},0,n_{0R}\right)
\end{array}}
&&
{\begin{array}{c}
\left(-n_{\Delta R},-n_{bR}-n_{0L},-n_{\Delta R}n_{0L},-n_{bR}n_{0L}\right)\\
\left(0,-n_{0L},-n_{\Delta R},-n_{bR}\right)
\end{array}}
\\
&
\left(c_{\Delta L},c_{bL},c_{\Delta R},c_{bR}\right)
\ar@{|->}`l[lu][lu]\ar@{|->}`r[ru][ru]
}\]
Also the map $d_1:E_1^{1,2n}\to E_1^{2,2n}$ is given by
\begin{align*}
&(a^1,a^2)\mapsto(n_{0L}+n_{0R},n_{0L}n_{0R})\\
&(n_{\Delta L},n_{bL},n_{0R})\mapsto(0,n_{0L},n_{0R})\\
&(n_{0L},n_{\Delta R},n_{bR})\mapsto(n_{0L},0,n_{0R})
\end{align*}
Now let
\begin{align*}
&K_{AL}={\rm Ker}(H^*(A_L)\to H^*(A_0))&
&K_{AR}={\rm Ker}(H^*(A_R)\to H^*(A_0))\\
&K_{NL}={\rm Ker}(H^*(N_L)\to H^*(N_0))&
&K_{NR}={\rm Ker}(H^*(N_R)\to H^*(N_0))
\end{align*}
Then
\begin{align*}
&H^*(A_L)\cong\Zz[a^1,a^2]\oplus K_{AL}\ ,\ 
H^*(A_R)\cong\Zz[a^1,a^2]\oplus K_{AR}\\
&H^*(C)\cong\Zz[n_{L},n_{R}]\oplus K_{NL}\oplus K_{NR}\oplus K_C\\
&H^*(N_L)\cong\Zz[n_{L},n_{R}]\oplus K_{NL}\ ,\ 
H^*(N_R)\cong\Zz[n_{L},n_{R}]\oplus K_{NR}
\end{align*}
Notice that $K_C\subset H^*(C)$ is the ideal generated by
$c_{\Delta L}c_{\Delta R}$. The restriction of the map
$H^*(A_L)\to H^*(N_L)$ to $K_{AL}$ induces a map
$s_L:K_{AL}\to K_{NL}$. Similarly we have a map $s_R:K_{AR}\to K_{NR}$.
Let also $s:\Zz[a^1,a^2]\to\Zz[n_L,n_R]$ be the map induced by the direct sum
map $BU(1)\times BU(1)\to BU(2)$.
Then the map $d_1:E_1^{0,2n}\to E_1^{1,2n}$
is given by
\begin{multline*}
d_1\left(a_L+k_{AL},a_R+k_{AR},x+k_{NL}+k_{NR}+k_C\right)=\\=
\left(-s(a_L)-s_L(k_{AL})+x+k_{NL},-s(a_R)-s_R(k_{AR})-x-k_{NR},a_L+a_R\right)
\end{multline*}
and the map $d_1:E_1^{1,2n}\to E_1^{2,2n}$ is given by
\[
d_1\left(x_L+k_{NL},x_R+k_{NR},a\right)=\left(x_L+x_R+s(a)\right)
\]
Now we can finish the proof:
\begin{enumerate}
\item We prove first that $H^{2n+1}(\Mm_2(X_2))=0$. We need to show
${\rm Ker}\,(d_1:E_1^{1,2n}\to E_1^{2,2n})\subset
{\rm Im}\,(d_1:E_1^{0,2n}\to E_1^{1,2n})$.
Let $\left(x_L+k_{NL},x_R+k_{NR},a\right)\in{\rm Ker}\,d_1$. Then
$x_L+x_R+s(a)=0$. It follows that
\[
d_1(a,0,-x_R+k_L-k_R)=(x_L+k_L,x_R+k_R,a)
\]
\item Now we will show that 
$H^{2n}\cong\Zz[a^1,a^2]\oplus K_{AL}\oplus K_{AR}\oplus K_C$ which 
conpletes the proof. We first define a map
$\Zz[a^1,a^2]\oplus K_{AL}\oplus K_{AR}\oplus K_C\to E_1^{0,2n}$ by
\[
\left(a,k_{AL},k_{AR},k_C\right)\mapsto
\left(a+k_{AL},-a+k_{AR},s(a)+s_L(k_{AL})-s_R(k_{AR})+k_C\right)
\]
We want to show this map is injective onto the kernel of $d_1$.
Injectivity is clear and a direct verification shows the image is contained
in the kernel of $d_1$. To show surjectivity let
$(a_L+k_{AL},a_R+k_{AR},x+k_{NL}+k_{NR}+k_C)\in{\rm Ker}\,d_1$.
Then 
\[
a_L=-a_R,\,k_{NL}=s_L(k_{AL}),\,k_{NR}=-s_R(k_{AR}),\,x=s(a_L)=-s(a_R)
\]
The result follows.
\end{enumerate}
\end{proof}

We are ready to prove the general case:

\begin{teor}
With notations as in theorem \ref{theo21} let
\begin{align*}
&K_i={\rm Ker}\left(\,
H^*(\pi_i^*\Mm_2(X_1))\to H^*(\pi_\emptyset^*\Mm_2(X_0))\,\right)\\
&K_{ij}={\rm Ker}\left(\,H^*(\pi_{ij}^*\Mm_2(X_2))\to
H^*(\pi_i^*\Mm_2(X_1))\oplus H^*(\pi_j^*\Mm_2(X_1))\,\right)
\end{align*}
Then, as modules over $\Zz$, we have an isomorphism
\begin{equation}\label{eq7main}
H^*\left(\,\Mm_2(X_q)\,\right)\cong 
H^*\left(\,\Mm_2(X_0)\,\right)\oplus \bigoplus_iK_i\oplus 
\bigoplus_{i<j}K_{ij}
\end{equation}
\end{teor}
\begin{proof}
We divide the proof into two steps:
\begin{enumerate}
\item We will use theorem
\ref{theo21} to build a spectral sequence converging to
the cohomology of $H^*(\Mm_2(X_q)$. Let $\Delta$ be the $q-1$ simplex.
Label its vertices by $v_i$,
$i=1,\ldots,q$, and the $e_{ij}$ be the middle point of the
edge joining $v_i$ and $v_j$. We define a filtration
$\Delta_0\subset\Delta_1\subset\Delta$ of $\Delta$ where
$\Delta_0=\bigcup_{i<j}e_{ij}$ and 
$\Delta_1$ is the 1-skeleton of $\Delta$. Write
$\Delta_1=\bigcup_i\Delta_{1i}$ where $\Delta_{1i}$ is the closure of
the connected component of $\Delta_1\setminus\Delta_0$ containing $v_i$.
Then we define
\[
M=\frac{
\bigcup_{i,j}\left(e_{ij}\times\pi_{ij}^*\Mm_2(X_2)\right)\cup
\bigcup_i\left(\Delta_{1i}\times\pi_i^*\Mm_2(X_1)\right)\cup
\left(\Delta\times\pi_\emptyset^*\Mm_2(X_0)\right)
}{\sim}
\]
where $\sim$ is induced by the inclusions
$e_{ij}\subset\Delta_{1i}\subset\Delta$ and 
$\pi_\emptyset^*\Mm_2(X_0)\subset\pi_i^*\Mm_2(X_1)\subset\pi_{ij}^*\Mm_2(X_2)$.
Then, the arguments in \cite{Seg68} can be applied to show that
$M$ is homotopically equivalent to $\Mm_2(X_q)$.
The filtration of $\Delta$ by $\Delta_0,\Delta_1$
induces a filtration $F_0\subset F_1\subset F_2=M$ of $M$ which leads
to a spectral sequence with
\begin{align*}
&E_1^{0,n}=H^n(F_0)\cong\bigoplus_{i<j}\left(
H^0(e_{ij})\otimes H^n(\pi_{ij}^*\Mm_2(X_2))\right)\\
&E_1^{1,n}=H^n(F_1,F_0)\cong\bigoplus_i\left(
H^1(\Delta_{1i},\partial\Delta_{1i})\otimes H^n(\pi_i^*\Mm_2(X_1))\right)\\
&E_1^{2,n}=H^n(F_2,F_1)\cong H^1(F_1)\otimes H^n(\Mm_2(X_0))
\end{align*}
\item The $d_1$ differential is induced by the inclusions
$\pi_\emptyset^*\Mm_2(X_0)\to\pi_i^*\Mm_2(X_1)\to\pi_{ij}^*\Mm_2(X_2)$.
We will use the sign convention
\begin{equation}\label{eq5nnnsignsq}
\xymatrix{
&H^*\left(\,\pi_i^*\Mm_2(X_1)\,\right)\ar[rd]^-{+}\\
H^*\left(\,\pi_{ij}^*\Mm_2(X_2)\,\right)\ar[ru]^-{+}\ar[rd]^-{-}&&
H^*\left(\,\pi_\emptyset^*\Mm_2(X_0)\,\right)\\
&H^*\left(\,\pi_j^*\Mm_2(X_1)\,\right)\ar[ru]^-{+}
}\end{equation}
Let
\begin{align*}
&K_i={\rm Ker}\left(\,
H^*(\pi_i^*\Mm_2(X_1))\to H^*(\pi_\emptyset^*\Mm_2(X_0))\,\right)\\
&K_{ij}={\rm Ker}\left(\,H^*(\pi_{ij}^*\Mm_2(X_2))\to
H^*(\pi_i^*\Mm_2(X_1))\oplus H^*(\pi_j^*\Mm_2(X_1))\,\right)
\end{align*}
Then, from theorem \ref{theo5nnnk=2} we have
\begin{align*}
&H^*(\pi_i^*\Mm_2(X_1))\cong H^*(\pi_\emptyset^*\Mm_2(X_0))\oplus K_i\\
&H^*(\pi_{ij}^*\Mm_2(X_2))\cong H^*(\pi_\emptyset^*\Mm_2(X_0))\oplus
K_i\oplus K_j\oplus K_{ij}
\end{align*}
Then the sequence of maps
$E_1^{0,n}\overset{d_1}{\rightarrow}
E_1^{1,n}\overset{d_1}{\rightarrow}E_1^{2,n}$
splits into three sequences
\begin{equation}\label{eq3seq}
\xymatrix @R=1em{
\bigoplus_{i<j}H^0(e_{ij})\otimes K_{ij}\ar[r]&0\ar[r]&0\\
\bigoplus_{i<j}H^0(e_{ij})\otimes(K_i\oplus K_j)\ar[r]&
\bigoplus_lH^1(\Delta_{1l},\partial\Delta_{1l})\otimes K_l\ar[r]&0\\
H^0(\Delta_0)\otimes K^n\ar[r]&
H^1(\Delta_1,\Delta_0)\otimes K^n\ar[r]&
H^1(\Delta_1)\otimes K^n
}\end{equation}
where $K^n$ stands for $H^n\left(\Mm_2(X_0)\right)$.
The bottom maps are easily analyzed using the
exact sequence
\[
0\to H^0(\Delta_1)\to H^0(\Delta_0)\to 
H^1(\Delta_1,\Delta_0)\to H^1(\Delta_1)\to0
\]
It follows that the map $d_1:E_1^{1,n}\to E_1^{2,n}$ is surjective. Since
$E_1^{r,n}=0$ for $r>2$ and $n$ even, this implies the spectral sequence
collapses and
\begin{gather*}
H^n(\Mm_2(X_q))=\frac{{\rm Ker}(d_1:E_1^{1,n}\to E_1^{2,n})}{
{\rm Im}(d_1:E_1^{0,n}\to E_1^{1,n})}\\
H^n(\Mm_2(X_q))={\rm Ker}(d_1:E_1^{0,n}\to E_1^{1,n})
\end{gather*}
Lets look more closely at the map
\begin{equation}\label{eq5nnnmiddle}
\bigoplus_{i<j}H^0(e_{ij})(K_i\oplus K_j)\to
\bigoplus_iH^1\left(\Delta_{1i},\partial\Delta_{1i}\right)\otimes K_i
\end{equation}
Observe that
\[
\bigoplus_{i<j}H^0(e_{ij})\otimes\left(\,
K_{i}\oplus K_{j}\,\right)=
\bigoplus_iH^0(\partial \Delta_{1i})\otimes K_i
\]
It follows that the map (\ref{eq5nnnmiddle}) can be easily analysed using the
exact sequence
\[
0\to H^0(\Delta_{1i})\to H^0(\partial\Delta_{1i})\to 
H^1(\Delta_{1i},\partial\Delta_{1i})\to0
\]
We gather together our conclusions:
\begin{enumerate}
\item The top sequence in (\ref{eq3seq}) contributes a term
\[
\bigoplus_{i<j}H^0(e_{ij})\otimes K_{ij}
\]
to $H^{2n}(\Mm_2(X_q))$.
\item The bottom sequence in (\ref{eq3seq}) does not contribute to 
$H^{2n+1}(\Mm_2(X_q))$ since it is exact in the middle.
\item The bottom sequence in (\ref{eq3seq}) contributes a term
\[
H^0(\Delta_1)\otimes H^*\left(\Mm_2(X_0)\right)
\]
to $H^{2n}(\Mm_2(X_q))$.
\item The map (\ref{eq5nnnmiddle}) is surjective hence it 
does not contribute to $H^{2n+1}(\Mm_2(X_q))$.
\item The map (\ref{eq5nnnmiddle}) contributes a term
\[
\bigoplus_iH^0(\Delta_{1i})\otimes K_i
\]
to $H^{2n}(\Mm_2(X_q))$.
\end{enumerate}
From (b) and (d) it follows that $H^{2n+1}(\Mm_2(X_q))=0$
and from (a), (c) and (e) equation (\ref{eq7main}) follows.
\end{enumerate}
\end{proof}

\bibliography{bi}

\end{document}